\documentclass[11pt,reqno,a4paper]{amsart}
\pdfoutput=1

\usepackage{iftex}
\ifPDFTeX
  \usepackage[utf8]{inputenc}
  \usepackage[T1]{fontenc}
  \usepackage{lmodern}
\else
  \usepackage[no-math]{fontspec}
\fi

\usepackage{amsmath,amsthm,amssymb,amsfonts,mathtools}
\usepackage{mathrsfs}
\usepackage[margin=1in]{geometry}
\usepackage{booktabs,array}
\usepackage[unicode=true,hidelinks,hypertexnames=false]{hyperref}
\hypersetup{pdfauthor={},pdftitle={Sixth and twelfth moments of Dirichlet L-functions to prime power moduli}}
\numberwithin{equation}{section}
\allowdisplaybreaks[2]
\theoremstyle{plain}
\newtheorem{theorem}{Theorem}[section]
\newtheorem{lemma}[theorem]{Lemma}
\newtheorem{proposition}[theorem]{Proposition}

\theoremstyle{remark}
\newtheorem{remark}[theorem]{Remark}
\date{}
\begin{document}

\title{Sixth and twelfth moments of Dirichlet $L$-functions\\ in prime power moduli}

\author{Xinchen Miao}
\address{Charles University, Faculty of Mathematics and Physics, Department of Algebra, Sokolovská 49/83, 186 75 Praha 8, Czech Republic}
\email{xinchen.miao@matfyz.cuni.cz, olivermiaoxinchen@gmail.com}

\author{Huimin Zhang}
\address{School of Mathematical Sciences, Zhejiang University \\
Hangzhou, Zhejiang 310058, China}
\email{huimin.zhang@zju.edu.cn}

\begin{abstract}
For a fixed odd prime $p$ and $q=p^n$, we establish
power-saving bounds for the sixth and twelfth moments
of central values of primitive Dirichlet $L$-functions,
with exponents $5/4-3/332$ and $2-51/261961$,
respectively. More generally, our methods yield
a power saving over the classical interpolation
bound for every moment of order $4<j\le12$.
Our approach begins with short second moments
over character cosets and preserves the dependence
on the character parameter. An exact linearization
of the $p$-adic stationary phases reduces the critical
large-value range to exponential sums with perturbed
square-root phases over prime-power rings. We establish
a large-value estimate for these sums by combining
congruence-spacing bounds, higher congruence energies,
and an average of quadratic Gauss sums over the dual
frequency. For the remaining high-value range, we
develop an incomplete-correlation estimate using
 $p$-adic exponent pairs. 
\end{abstract}



\maketitle

\makeatletter
\def\ps@headings{\ps@empty
  \def\@evenhead{\normalfont\scriptsize\rlap{\thepage}\hfil MOMENTS OF DIRICHLET $L$-FUNCTIONS\hfil}
  \def\@oddhead{\normalfont\scriptsize\hfil MOMENTS OF DIRICHLET $L$-FUNCTIONS\hfil\llap{\thepage}}
  \let\@mkboth\markboth}
\makeatother
\pagestyle{headings}

\section{Introduction}\label{sec:intro}

Moment estimates for the Riemann zeta function and Dirichlet $L$-functions
play a fundamental role in analytic number theory.
For $\ell\ge0$, define
$$
M_\ell(T):=\int_0^T|\zeta(1/2+it)|^{\ell}\,dt.
$$
The Lindel\"of hypothesis is equivalent to the
assertion that $M_\ell(T)\ll_{\ell,\varepsilon}T^{1+\varepsilon}$
for every fixed $\ell\ge0$ and every $\varepsilon>0$.
Although the  Lindel\"of-strength bound is known unconditionally only for $0\le \ell \le4$, upper bounds for higher moments provide strong control over the large-value distribution of the Riemann zeta function. Obtaining improvements over the classical
interpolation bounds remains an important problem. They are also closely related to subconvexity problems.

The classical fourth moment estimate
$
M_4(T)\ll_\varepsilon T^{1+\varepsilon}
$
was established by Hardy and Littlewood \cite{HL} and subsequently Ingham
\cite{Ing} proved the asymptotic formula.
A fundamental higher-moment
result is Heath-Brown's twelfth moment bound
\cite{HB}
$
M_{12}(T)\ll_\varepsilon T^{2+\varepsilon}.
$
Interpolating between these two estimates by
H\"older's inequality gives
$M_\ell(T)\ll_\varepsilon
T^{1+(\ell-4)/8+\varepsilon}$ 
for $ 4\le \ell\le12 $.
In particular, this yields the classical
sixth moment bound
$
M_6(T)\ll_\varepsilon T^{5/4+\varepsilon}.
$
Recent progress has improved these exponents.
In particular, de Faveri and Pandey \cite{DFP}
obtained an improved sixth moment estimate
$
M_6(T)\ll_\varepsilon
T^{5/4-1/60+\varepsilon},
$
and, more generally, improved the upper bounds
for $M_\ell(T)$ throughout the range $4<\ell\le12$.
Their approach is based on 
estimates for large values of the Riemann zeta function and a large-value estimate
for exponential sums with square-root phases and certain perturbations.

Motivated by these developments, we investigate
corresponding improvements for Dirichlet
$L$-functions in the depth aspect.

Let $p$ be a fixed odd prime and let $q=p^n$,
$n\ge1$. For $j\ge1$, define
\begin{equation}\label{intro:moments}
M_j(q)=\sum_{\chi\bmod q}^{*}|L(1/2,\chi)|^j,
\end{equation}
where the star denotes primitive Dirichlet
characters of conductor exactly $q$. The depth
aspect, in which $p$ is fixed and $n$ tends to
infinity, provides a finite-place counterpart
to the $t$-aspect of the Riemann zeta
function.

The study of moments for Dirichlet $L$-functions in depth aspect was advanced
by Mili\'cevi\'c and White \cite{MW}, who proved
the prime-power analogue of Heath-Brown's
twelfth moment bound for the Riemann zeta
function \cite{HB}. Together with the fourth moment
estimate for Dirichlet $L$-functions
\cite{Sound}, their results give
\begin{equation}\label{intro:classical}
M_4(q)\ll_{p,\varepsilon}q^{1+\varepsilon},
\qquad
M_{12}(q)\ll_{p,\varepsilon}q^{2+\varepsilon},
\qquad
M_6(q)\ll_{p,\varepsilon}q^{5/4+\varepsilon}.
\end{equation}
The last estimate follows by H\"older's inequality.
The corresponding moment problems have also been
studied for other families and moduli, including
smooth square-free moduli and related fourth moment
problems; see, for example,
\cite{HB1,HB2,Irv,Moto,Nun,PY1,PY2}. 

Pointwise estimates provide a complementary
approach to the study of central values.
In the $t$-aspect, Bourgain \cite{Bou} obtained a
sub-Weyl bound for the Riemann zeta function using
decoupling and exponential sum estimates.
In the depth aspect, the pointwise bound of
Mili\'cevi\'c \cite{Mil} gives, for fixed odd $p$,
\begin{equation}\label{eq:inputs}
|L(1/2,\chi)|
\ll_{p,\varepsilon}q^{27/164+\varepsilon}.
\end{equation}
This estimate will be used to provide an upper
cutoff for the central values in the proof of the twelfth moment bound.

Our aim is to improve the sixth and
twelfth moment exponents in the depth aspect.
We obtain the following result.

\begin{theorem}[Sixth and twelfth moment bounds]
\label{thm:main}
Let $p$ be a fixed odd prime and let $q=p^n$, $n\ge1$. For every
$\varepsilon>0$, we have
\begin{equation}\label{eq:main}
\begin{aligned}
M_6(q)&\ll_{p,\varepsilon}
q^{5/4-3/332+\varepsilon},\\
M_{12}(q)&\ll_{p,\varepsilon}
q^{2-51/261961+\varepsilon}.
\end{aligned}
\end{equation}
\end{theorem}
Thus, the savings over the classical exponents
$5/4$ and $2$ are $3/332$ and $51/261961$,
respectively. In particular,
$5/4-3/332=103/83$.

Interpolation between the fourth moment estimate
and Theorem~\ref{thm:main} gives the following
general moment bound.

\begin{theorem}[General moment bound]\label{holder}
For each $4<j\le12$, there exists $\delta_j>0$
such that, for every $\varepsilon>0$ and
$q=p^n$, with $p$ a fixed odd prime and $n\ge1$,
\begin{equation}\label{eq:holder}
M_j(q)\ll_{p,\varepsilon}
q^{1+(j-4)/8-\delta_j+\varepsilon}.
\end{equation}
\end{theorem}

\begin{proof}
For $j=6$ and $j=12$, Theorem~\ref{thm:main}
gives the explicit savings
$\delta_6=3/332$ and
$\delta_{12}=51/261961$, respectively.
For the remaining values of $j$, the existence
of a positive $\delta_j$ follows by interpolation
between the fourth moment estimate in \cite{Sound} and the new twelfth moment bound in Theorem \ref{thm:main}.
\end{proof}

\begin{remark}
The methods developed in this paper yield
further improvements beyond the range of
Theorem~\ref{holder}. In particular, improved
moment bounds can also be obtained for
$12<j<12+\nu_0$, where $\nu_0>0$ is small
(for example, one may take $\nu_0\approx0.024$).
We do not pursue the optimization of either
the savings or this range here.
For example, we can improve $\delta_{12}$ by carefully picking the positive integer $k$ (may pick $k=20$ instead of $k=19$ in the paper) and a suitable exponent pair.
For $4<j<12$ with
$j\ne6$, applying the methods of
Sections~\ref{sec:finite} and
\ref{sec:counting-large-values}, rather than
direct interpolation, also gives better savings.
For positive integer moments, some choices of
the energy order $k$ and the resulting savings
are as follows:
\begin{align*}
&k=3,\qquad
\delta_5=\frac3{664},\quad
\delta_6=\frac3{332},\quad
\delta_7=\frac9{664},\quad
\delta_8=\frac3{166};\\
&k=3,4,\qquad
\delta_9=\frac{59}{3000};\\
&k=6,7,\qquad
\delta_{10}=\frac{109}{7308};\\
&k=13,14,\qquad
\delta_{11}=\frac{601}{75320}.
\end{align*}
Here the admissible values of $k$ are specified
in Theorem~\ref{thm:finite}.
\end{remark}

\subsection{The critical range and the short-family method}
\label{intro:MW}

We first recall the short-family method of
Mili\'cevi\'c and White \cite{MW}, which provides
the basic reduction for the 
depth aspect. Let $R(q,V)$ denote the number
of primitive characters modulo $q$ satisfying
$|L(1/2,\chi)|\ge V$. The fourth moment estimate
in \cite{Sound} and the twelfth moment bound
in \cite{MW} imply
\begin{equation}\label{intro:classical-count}
R(q,V)\ll_{p,\varepsilon}q^\varepsilon
\min\{qV^{-4},q^2V^{-12}\}.
\end{equation}
The two bounds coincide at $V=q^{1/8}$.
A dyadic block
$V\le |L(1/2,\chi)|<2V$ contributes
$O(V^6R(q,V))$ to the sixth moment.
Consequently, the transition at $V=q^{1/8}$
gives the classical exponent $5/4$.
For the twelfth moment, the second term in
\eqref{intro:classical-count} gives
$V^{12}R(q,V)\ll q^{2+\varepsilon}$.
Thus an improvement near the critical range
$q^{1/8}$ is necessary for the sixth moment,
whereas a saving for the twelfth moment also
requires a further improvement for larger values.

For a divisor $d\mid q$, let $H_d$ denote the
full group of characters modulo $d$, inflated
to modulus $q$. Following \cite{MW}, we consider
the short second moment
\begin{equation}\label{intro:short}
S_2(\chi;q_1)
=
\sum_{\psi\in H_{q_1}}
|L(1/2,\chi\psi)|^2,
\qquad p\le q_1<q.
\end{equation}
The sum is over all lower-conductor twists,
not only primitive ones. When $\chi$ has
conductor $q$, every character in the coset
$\chi H_{q_1}$ remains primitive. We choose
$q_1$ slightly smaller than $V^2$, so that a
coset containing $F$ characters with central
values of size at least $V$ contributes at
least $FV^2$ to the short second moment.

The approximate functional equation, character
orthogonality and Poisson summation transform
this short moment into sums involving the
complete kernels $K_A(m;Q_h)$ defined in
Section~\ref{sec:prelim}, with
\begin{equation}\label{intro:scales}
D_h=q_1h,\qquad
Q_h=q/D_h,\qquad
N_h=q/D_h^2,\qquad h=p^j.
\end{equation}
This reduction, including its conductor descent,
is the starting point supplied by
\cite[Proposition~1]{MW}. The important feature
for our purposes is that the resulting
coefficients can be chosen independently of
the center character.

Mili\'cevi\'c and White estimate the resulting
sums through pairwise correlations and
completion. For $q_1\mid q_2<q$ and
$\Psi\subseteq H_{q_2}$, their aggregate estimate
\cite[Proposition~3]{MW} is
\begin{equation}\label{intro:MW-aggregate}
\begin{aligned}
\sum_{\psi\in\Psi}S_2(\chi\psi;q_1)
\ll_{p,\varepsilon}q^\varepsilon
\bigl(
(q_1+q_1^{1/4}q_2^{1/4})|\Psi|
+q^{1/2}|\Psi|^{1/2}
\bigr).
\end{aligned}
\end{equation}
With $q_1$ of size $V^2$, absorption of the
first term leads to $q_2\le q_1^3$, up to
factors of $q^\varepsilon$. This recovers
the second bound in
\eqref{intro:classical-count}. Our argument
uses the same short-moment reduction but
treats the critical and higher ranges by
different estimates.

Sections~\ref{intro:finite} and \ref{intro:tail} outline the proof of Theorem~\ref{thm:main}.

\subsection{A finite-ring large-value argument}
\label{intro:finite}

We first improve the large-value estimate near
the critical threshold $V=q^{1/8}$. Choose a
power $q_1$ of $p$ slightly smaller than $V^2$.
If an $H_{q_1}$-coset contains $F$ characters
with $|L(1/2,\chi)|\ge V$, then its short
second moment is at least $FV^2$. By
Lemma~\ref{lem:short}, after fixing a coefficient
system, a conductor-descent level $h$, a
stationary sign and a square-root branch, we
obtain
\begin{equation}\label{large-value}
\left|
\sum_{m\in I}a_mK_A^\pm(m;Q_h)
\right|
\gg q^{-\varepsilon}FQ_h^{1/2},
\end{equation}
where $D_h=q_1h$, $Q_h=q/D_h$ and
$N_h=q/D_h^2$. Here $I$ is an ordinary integer
interval of length $O(q^\varepsilon N_h)$,
and the coefficients $a_m$ are independent of
the center character.

Instead of completing pairwise correlations
as in \cite{MW}, we preserve the dependence on
the Postnikov parameter $A$. Fix an
$H_{q_2}$-coset with $q_2=D_hR<q$ and let
$\Delta=q/q_2$, so that $Q_h=\Delta R$.
On a fixed square-root branch, write
$A=c(B_0+\Delta w)^2$ and $m=cy^2$, where
$c,B_0$ are fixed units, $w$ runs modulo $R$,
and $y=y(m)$ is independent of $w$.
Proposition~\ref{prop:phase} gives the exact
congruence
$$
\Phi_{A,D_h}^\pm(m)
\equiv C^\pm(m)
\pm2c\Delta w(y+\eta y^3)
\pmod{Q_h},
$$
where $C^\pm(m)$ is independent of $w$ and
$\eta=-D_h^2/(24B_0^2)$. The normalized
quadratic Gauss factor is also independent of
$w$ on the fixed branch. Absorbing these
fixed factors into the coefficients and
rescaling $w$ by the unit $\pm2c$, we obtain
\begin{equation}\label{intro:polynomial}
\mathcal D(w)
=
\sum_{m\in I}a_m
e\left(\frac{w(r(m)+\eta r(m)^3)}R\right),
\qquad
r(m)^2=m/c,\qquad
\eta^2\equiv0\pmod R.
\end{equation}
The coefficients $a_m$ may depend on the
chosen coset and square-root branch, but are
fixed as $w$ varies. After normalization,
we may assume $|a_m|\le1$. The cubic
perturbation cannot in general be discarded.
The nilpotence condition $\eta^2\equiv0\pmod R$
makes the map $x\mapsto x+\eta x^3$ invertible
modulo $R$, which is essential for the subsequent
congruence analysis.

Let $\mathcal W\subset\mathbb Z/R\mathbb Z$
be a set of parameters satisfying
$|\mathcal D(w)|\ge Z$. Theorem~\ref{thm:finite}
bounds $|\mathcal W|$ by combining two estimates.
First, the congruence-spacing lemma for the
frequencies $r(m)+\eta r(m)^3$, together with
the large-value condition, gives the higher-energy
estimate
\begin{equation*}
\mathcal E_{2k}(D_0)
\ll_{p,k,\varepsilon}R^\varepsilon
\left(
\left(\frac NZ\right)^{4k-4}\mathcal E_2(D_0)
+\frac{ND_0}{Z^2}|\mathcal W|^{2k-1}
\right).
\end{equation*}
Here $\mathcal E_{2k}(D_0)$ is defined in
\eqref{eq:energy-def}, and
$\mathcal E_2(D_0)\le D_0|\mathcal W|$.
Second, Poisson summation in the original
summation variable produces, after conductor
descent, the quadratic--quartic phase
$$
\Psi_b(z)=\frac{z^2}{4cb}
+\frac{\eta z^4}{8c^3b^3}.
$$
Averaging the corresponding complete sums
over the dual frequency $b$ yields the weight
$D_0^{-1/2}$ in the geometric estimate of
Lemma~\ref{lem:geometric}. Combining this
geometric estimate with the higher-energy
bound proves Theorem~\ref{thm:finite}.

The real-variable model for this argument is
provided by de Faveri and Pandey
\cite[Proposition~4.1 and Sections~5--6]{DFP},
who obtain large-value estimates for exponential
sums with perturbed square-root phases by
combining higher additive energies with
pointwise coefficient bounds. The use of
large-value sets and coefficient-sensitive
energy estimates is also related to the
large-value methods of Guth and Maynard
\cite[Section~11]{GM}. The present finite-ring
setting requires separate arguments for
congruence spacing, exact local phase
linearization and the treatment of non-primitive
Poisson frequencies. In particular, the
real-variable estimates of \cite{DFP} cannot
be applied directly to the prime-power setting.

The improvement is already visible at the
classical threshold. Suppressing arbitrarily
small powers of $q$, the principal level
$h=1$ has
\begin{equation}\label{intro:critical-scales}
V=q^{1/8},\quad q_1=q^{1/4},\quad
N=q^{1/2},\quad Q=q^{3/4},\quad
R=q^{7/12},\quad Z=q^{3/8},\quad
\Upsilon=q^{1/12}.
\end{equation}
Here $\Upsilon=R/(R,\eta)$, where
$\eta=-D_h^2/(24B_0^2)$ as above.
Since $h=1$ and $B_0$ is a $p$-adic unit,
we have $(R,\eta)\asymp_p q^{1/2}$,
up to the suppressed $q^\varepsilon$ losses.
Theorem~\ref{thm:finite} with $k=3$ gives
at most $q^{5/16+\varepsilon}$ large parameters
in each fixed $H_{q_2}$-coset. Since there are
$O_p(q/(q_1R))$ such cosets, we obtain
\begin{equation}\label{intro:critical-gain}
R(q,q^{1/8})\ll_{p,\varepsilon}q^{23/48+\varepsilon},
\qquad
\frac{23}{48}=\frac12-\frac1{48}.
\end{equation}
This improvement comes from the distribution
of large exponential sums, rather than from
a stronger pointwise bound for individual
central values.

Section~\ref{sec:counting-large-values} extends
the estimate to neighboring parameter ranges,
takes into account the exact $h$-fold
multiplicity of the parameter fibers, and shows
that the contribution from every
conductor-descent level is bounded by that
from the level $h=1$. The choice $k=3$ yields
the critical-range estimate used for the sixth
moment, while $k=19$ extends the saving to
the upper transition range needed for the
twelfth moment.

\subsection{The higher range and the two moment bounds}
\label{intro:tail}

The finite-ring large-value argument controls
the critical range but does not cover the
entire range needed for the twelfth moment.
For larger values, we return to the short-moment
aggregation in \eqref{intro:MW-aggregate}
and improve its correlation estimate before
completion. Let $A,A'$ be the Postnikov
parameters of two characters \cite{Pos} in the same
$H_{q_2}$-coset, and let
\begin{equation}\label{intro:correlation-scales}
P_h=\frac{Q_h}{(Q_h,A-A')},\qquad
R_0=q_2/q_1,\qquad N_0=q/q_1^2,
\end{equation}
where $P_h$ is the effective conductor of the
phase difference. Lemma~\ref{lem:corr} bounds
the incomplete same-sign correlation over an
interval of length $M$ by
$q^\varepsilon P_h^\kappa M^{\lambda-\kappa}$,
for any fixed admissible $p$-adic exponent datum
with exponent pair $(\kappa,\lambda)$ and
sufficiently large $P_h$. The proof normalizes
the difference of the stationary phases and
verifies the full exponent-data hypotheses
of \cite{Mil}. The remaining small-conductor
correlations are controlled by the trivial
bound and the distribution of nearby character
parameters.

Applying Cauchy--Schwarz and summing over the
conductor-descent levels gives the improved
aggregation estimate
\begin{equation}\label{intro:new-aggregate}
\sum_{\psi\in\Psi}S_2(\chi\psi;q_1)
\ll  q^\varepsilon \left(
\left(
q_1+q_1^{1/2}R_0^{\kappa/2}
N_0^{(\lambda-\kappa)/2}
\right)|\Psi|
+q^{1/2}|\Psi|^{1/2} \right).
\end{equation}
The classical pair $(\kappa,\lambda)=(1/2,1/2)$
recovers \eqref{intro:MW-aggregate}. More
generally, the additional term proportional
to $|\Psi|$ is absorbed by $q_1|\Psi|$
whenever
$
R_0^\kappa N_0^{\lambda-\kappa}\le q_1.
$

To exploit this condition, write $V=q^\sigma$
and choose $q_1\asymp_p V^2$, with a small
power saving for absorption. We seek to take
$q_2$ of size $q_1^3q^\delta$, which gives
$R_0\asymp_p q^{4\sigma+\delta}$ and
$N_0\asymp_p q^{1-4\sigma}$, up to the
auxiliary small powers of $q$. The fixed
exponent pair
$(\kappa,\lambda)=(2346/4789,2400/4789)$
used in Section~\ref{tail:count} reduces the
absorption condition to
$
\delta\le\frac{205\sigma-27}{1173}
$
at the exponent level. Thus a positive saving
is available in the required higher range.

Applying \eqref{intro:new-aggregate} to the
short cosets containing large central values
gives
$$
R(q,V)\ll
\frac{q^{2+\varepsilon}}{q_1q_2V^4}.
$$
Since $q_1\asymp_p V^2$ and
$q_2\asymp_p q_1^3q^\delta$, this yields
$$
R(q,V)\ll_{p,\varepsilon}
q^{2-\delta+\varepsilon}V^{-12}.
$$
This is the large-value improvement over
\cite{MW}, complementary to the critical-range
energy argument.

For clarity, the bounds used in the final
summation are
\begin{equation}\label{intro:distribution}
\begin{aligned}
R(q,q^\sigma)
&\ll_{p,\varepsilon}
q^{(28-155\sigma)/18+\varepsilon},
&&\frac{10}{83}\le\sigma\le\frac{13}{100},\\
R(q,q^\sigma)
&\ll_{p,\varepsilon}
q^{(1404-8315\sigma)/738+\varepsilon},
&&\frac18\le\sigma\le\sigma_*,\\
R(q,V)
&\ll_{p,\varepsilon}
q^{2-\delta_*+\varepsilon}V^{-12},
&&q^{\sigma_*}\le V\le q^\beta,
\end{aligned}
\end{equation}
with $$
\sigma_*=\frac{34794}{261961},\qquad
\delta_*=\frac{51}{261961},\qquad
\beta=\frac16-\frac1{1000}.$$
These estimates are established in
Propositions~\ref{prop:critical-range} and
\ref{prop:tail}. The first large-value estimate
meets the classical bound
$R(q,q^\sigma)\ll q^{1-4\sigma+\varepsilon}$,
which follows from the fourth moment estimate,
at $\sigma=10/83$. This yields
$M_6(q)\ll q^{103/83+\varepsilon}$; its
higher-value tail is controlled by the
previously known twelfth moment \eqref{eq:inputs}.

For the improved twelfth moment, the
critical-range estimate with $k=19$ gives
the saving $(72-541\sigma)/738$, while the
exponent-pair method permits the saving
$(205\sigma-27)/1173$. Balancing these
savings determines the transition point
$\sigma_*$ and the common saving $\delta_*$.
The classical fourth moment controls the
smaller values, and the pointwise sub-Weyl
bound of \cite{Mil} excludes values beyond
$q^\beta$. Hence we have
$$
M_{12}(q)\ll_{p,\varepsilon}
q^{2-\delta_*+\varepsilon}.
$$
The sixth and twelfth moment deductions use
separate ranges and estimates; in particular,
the improved twelfth moment bound is not
used in its own proof.

In summary, Sections~\ref{sec:prelim}--\ref{sec:short}
establish the arithmetic reduction. Section
\ref{sec:finite} proves the finite-ring theorem,
while Sections~\ref{sec:counting-large-values}--\ref{sec:high-values} establish
\eqref{intro:distribution}. Section
\ref{sec:moments} completes the moment
summations.

\textbf{Notation}
Throughout, $p$ denotes a fixed odd prime and
$q=p^n$. We write
$e(x)=\exp(2\pi ix)$ and $v_p(0)=+\infty$.
For $x\in\mathbb Z_p$ and $S=p^\nu$, the expression
$e(x/S)$ is understood using any integer
representative of $x\bmod S$. We also write
$(x,S)=p^{\min\{v_p(x),\nu\}}.$ $\log_p$ denotes the $p$-adic
logarithm on $1+p\mathbb Z_p$, whereas $\log$
denotes the natural logarithm of a positive
real number. In particular,
$\log_q x=\log x/\log q$ denotes the real
logarithm to base $q$, for $x>0$.
For an integrable function $f$ on $\mathbb R$,
we use the Fourier transform convention
$
\widehat f(\xi)
=
\int_{\mathbb R}f(x)e(-x\xi)\,dx.
$
A superscript $*$ on a character sum or maximum
over characters modulo $q$ restricts to primitive
characters of conductor exactly $q$. On a sum
over residue classes modulo a prime power, a
superscript $*$ restricts to units.
Whenever a square root is used, its branch is
specified by the corresponding congruence
condition. Inverses are taken in the relevant
$p$-adic unit group.
The letters $D_h$ and $D_0$ denote, respectively,
an arithmetic conductor parameter and an energy
scale. We use $\mathfrak d$ for a fixed $p$-adic
exponent datum, while the ordinary letter $d$
denotes a locally defined divisor, valuation,
or exponent. 
Finally, $\varepsilon>0$ denotes an arbitrarily
small positive number, not necessarily the same
at each occurrence.

{\textbf{Acknowledgments}} 
X. Miao is supported by the Czech 
Science Foundation GAČR, grant 26-20514S and the Charles University programme PRIMUS/25/SCI/008.
H. Zhang is in part supported by the China Postdoctoral Science Foundation (No. 2026M793404) and the National Key R$\&$D Program of China (No. 2021YFA1000700). 

{\textbf{AI Usage}}
We used ChatGPT to assist on the computation and optimatzation in some parts of Section \ref{sec:counting-large-values} and \ref{sec:high-values}, especially the choice of certain suitable exponent pair in Section \ref{tail:count}. ChatGPT were also used to improve the exposition in parts of the manuscript and check the manuscript. We take full responsibility for the mathematical arguments, results, and final text.

\section{Preliminaries}\label{sec:prelim}

\subsection{Approximate functional equation}
We use a dyadic form of the approximate
functional equation for primitive Dirichlet
$L$-functions, following
\cite[Lemma~1]{MW}. It follows from the
standard approximate functional equation
\cite[Theorem~5.3]{IwaniecKowalski2004} by a smooth dyadic
partition of unity.
\begin{lemma}\label{AFE2}
Let $\chi$ be a primitive Dirichlet character modulo $q$. Then,
$$
    \big|L\big(\tfrac{1}{2},\chi\big)\big|^2 \ll \log q \sum_{\substack{N \leqslant q^{1/2 + \varepsilon} \\ N \textnormal{ dyadic}}} \left| \frac{1}{\sqrt{N}} \sum_n \chi(n) V_N(n) \right|^2+ q^{-100},
$$
where $V_N$ is a smooth function depending only on $N$ and $q$, whose support is contained in $[N/2,2N]$ and whose derivatives satisfy $V_N^{(j)} \ll_j N^{-j}$ for every $j\in\mathbb{N}$.
\end{lemma}

\subsection{Character parameters and stationary branches}\label{prelim:characters}

For odd $p$, the $p$-adic logarithm identifies $1+p\mathbb Z_p$ with $p\mathbb Z_p$. By \cite[Lemma~2]{MW}, a primitive character $\chi\bmod q$ has a unit Postnikov parameter $A$, unique modulo $q/p$, such that

\begin{equation}
\label{eq:postnikov}
\chi(1+px)=e\!\left(
\frac{A\log_p(1+px)}q
\right),
\qquad x\in\mathbb Z_p,\quad
A\in\mathbb Z_p^\times.
\end{equation}

For a $p$-power $d_0$ with $p\le d_0<q$,
recall that $H_{d_0}$ is the group of all
characters modulo $d_0$, inflated to modulus
$q$. If $\psi\in H_{d_0}$, its wild parameter
can be written as $\mathfrak b_\psi$ modulo
$d_0/p$, and
\begin{equation}\label{prelim:parameter-increments}
A_{\chi\psi}\equiv
A_\chi+\frac q{d_0}\mathfrak b_\psi
\pmod{q/p}.
\end{equation}
Each value of $\mathfrak b_\psi$ has exactly
$p-1$ lifts in $H_{d_0}$, corresponding to
the $p-1$ possible tame characters modulo $p$.
Since $d_0<q$, every $\psi\in H_{d_0}$ is
trivial on $1+p^{n-1}\mathbb Z_p$, whereas
$\chi$ is nontrivial there.
Hence, every
character in $\chi H_{d_0}$ is primitive.

For $p\le D\mid d_0$, let $Q=q/D$ and consider
the map
$
\psi\longmapsto A_{\chi\psi}\pmod Q.
$
By the preceding congruence, this map depends
on $\mathfrak b_\psi$ only modulo $d_0/D$.
Hence its image has size $d_0/D$. Each element
of the image has $D/p$ possible wild parameters,
each with $p-1$ tame lifts. Thus every fiber
has
$
\frac{D}{p}(p-1)=\varphi(D)
$
elements. These fiber sizes will be used to
count short cosets and pairs of nearby
parameters. Moreover, multiplication by an
element of $H_{q_1}$ permutes $H_{q_1}$, so
$S_2(\chi;q_1)$ is constant on the coset
$\chi H_{q_1}$.

For a $p$-power $D$ with $p\le D<q$, let
$Q=q/D$. The normalized complete kernel is
\begin{equation}\label{eq:kernel}
K_A(m;Q)=Q^{-1/2}\sum_{x\bmod Q}^{*}
 e\!\left(\frac{\frac A D\log_p(1+Dx)+\frac m x}{Q}\right).
\end{equation}
It depends only on $A\bmod Q$, since
$D^{-1}\log_p(1+Dx)\in\mathbb Z_p$.
For $Q\ge p^2$ and $p\nmid m$,
\cite[Lemma~6]{MW} shows that $K_A(m;Q)=0$
unless $m/A$ is a square modulo $p$.
Otherwise, the two stationary branches give
\begin{equation}\label{eq:branches}
\begin{gathered}
K_A(m;Q)=K_A^+(m;Q)+K_A^-(m;Q),\\
K_A^\pm(m;Q)
=
\omega_A^\pm(m;Q)
e\left(\frac{\Phi_{A,D}^\pm(m)}Q\right),
\qquad |\omega_A^\pm(m;Q)|=1,\\
x_\pm=
\frac{Dm\pm\sqrt{D^2m^2+4Am}}{2A},
\qquad
\Phi_{A,D}^\pm(m)
=
\frac AD\log_p(1+Dx_\pm)+\frac m{x_\pm}.
\end{gathered}
\end{equation}
Here the square root is chosen on a fixed
$p$-adic branch, so that $x_+$ and $x_-$
reduce to the two roots of
$x^2\equiv m/A\pmod p$, respectively.
The signs thus label the two
stationary branches. By
\cite[Lemma~4 and Remark~1]{MW}, the normalized
quadratic Gauss factor at a stationary point
depends only on its residue class modulo $p$,
on $A\bmod p$, and on the parity of $v_p(Q)$.
In particular, $|K_A(m;Q)|\le2$.
We write $K_\chi$ for $K_A$ when $A$ is
the parameter of $\chi$.

\subsection{Exponential sums over prime power rings}
We record three elementary estimates for exponential
sums over prime power rings, which will be used in
the subsequent spacing and Poisson summation arguments.
Throughout this subsection, let $p$ be an odd prime
and $S=p^\nu$ with $\nu\ge1$.
\begin{lemma}[Complete quadratic Gauss sums]\label{Complete quadratic Gauss sums}
For $\alpha,\beta\in\mathbb Z_p$, we define $g=(\alpha,S)$, with
$(0,S)=S$. Then
\begin{equation}\label{frequency:quadratic-formula}
\left|\sum_{x\bmod S}e\!\left(\frac{\alpha x^2+\beta x}{S}\right)\right|
=\begin{cases}
\sqrt{Sg},&g\mid\beta,\\
0,&g\nmid\beta.
\end{cases}
\end{equation}
\end{lemma}

\begin{proof}
We define $Q=S/g$. Since $S\mid\alpha Q$, translation by $Q$
gives
\begin{equation}
\sum_{x\bmod S}e\!\left(\frac{\alpha x^2+\beta x}{S}\right)
=e\!\left(\frac{\beta}{g}\right)
\sum_{x\bmod S}e\!\left(\frac{\alpha x^2+\beta x}{S}\right).
\end{equation}
Thus the sum vanishes if $g\nmid\beta$.

Suppose that $g\mid\beta$ and write
$\alpha=g\alpha_0$, $\beta=g\beta_0$. Conductor descent yields
\begin{equation}
\sum_{x\bmod S}e\!\left(\frac{\alpha x^2+\beta x}{S}\right)
=g\sum_{x\bmod Q}e\!\left(\frac{\alpha_0x^2+\beta_0x}{Q}\right).
\end{equation}
If $Q=1$, the sum equals $S=\sqrt{Sg}$.
If $Q>1$, then $(\alpha_0,Q)=1$; since $p$ is odd,
completing the square gives
\begin{equation}
\left|\sum_{x\bmod Q}e\!\left(\frac{\alpha_0x^2+\beta_0x}{Q}\right)\right|
=\left|\sum_{x\bmod Q}e\!\left(\frac{\alpha_0x^2}{Q}\right)\right|.
\end{equation}
The substitution $(r,s)=(x-y,x+y)$ is bijective modulo $Q$.
Therefore, by orthogonality,
\begin{equation}
\begin{aligned}
\left|\sum_{x\bmod Q}e\!\left(\frac{\alpha_0x^2}{Q}\right)\right|^2
&=\sum_{x,y\bmod Q}e\!\left(\frac{\alpha_0(x^2-y^2)}{Q}\right)\\
&=\sum_{r,s\bmod Q}e\!\left(\frac{\alpha_0rs}{Q}\right)=Q.
\end{aligned}
\end{equation}
Consequently, the original sum has absolute value
$g\sqrt Q=\sqrt{Sg}$.
\end{proof}

\begin{lemma}[Non-stationary cancellation]\label{lem:nonstationary-cancellation}
Let $u\in\mathbb Z_p^\times$, $h\in\mathbb Z_p[x]$, and
$f(x)=ux+ph(x)$. Then, for every $\nu\ge1$,
\begin{equation}\label{prelim:nonstationary-sum}
\sum_{x\bmod p^\nu}e\!\left(\frac{f(x)}{p^\nu}\right)=0.
\end{equation}
\end{lemma}
\begin{proof}
Write $x=y+p^{\nu-1}t$, with $y\bmod p^{\nu-1}$ and
$t\bmod p$. Since $h\in\mathbb Z_p[x]$,
\begin{equation}
f(y+p^{\nu-1}t)\equiv f(y)+up^{\nu-1}t\pmod{p^\nu}.
\end{equation}
Therefore, we have
\begin{equation}
\begin{aligned}
\sum_{x\bmod p^\nu}e\!\left(\frac{f(x)}{p^\nu}\right)
&=\sum_{y\bmod p^{\nu-1}}e\!\left(\frac{f(y)}{p^\nu}\right)
\sum_{t\bmod p}e\!\left(\frac{ut}{p}\right)=0,
\end{aligned}
\end{equation}
as $p\nmid u$.
\end{proof}

\begin{lemma}[Quadratic stationary phase \cite{GK}]\label{lem:gauss}
Let
$f(t)=a_2t^2+a_1t+a_0+ph(t)$, where
$a_0,a_1\in\mathbb Z_p$, $a_2\in\mathbb Z_p^\times$, and
$h\in\mathbb Z_p[t]$. For $\alpha,\beta\in\mathbb Z_p$ with
$\alpha\not\equiv0\pmod S$, we let $g=(\alpha,S)$. Then we have 
\begin{equation}\label{eq:gauss-basic}
\begin{gathered}
\left|\sum_{t\bmod S}
e\!\left(\frac{\alpha f(t)+\beta t}{S}\right)\right|
\le\sqrt{Sg},\\
g\nmid\beta\quad\Longrightarrow\quad
\sum_{t\bmod S}e\!\left(\frac{\alpha f(t)+\beta t}{S}\right)=0.
\end{gathered}
\end{equation}
\end{lemma}
\begin{proof}
We let $Q=S/g$ and write $t=y+Qv$, with $y\bmod Q$ and
$v\bmod g$. Since $f(y+Qv)\equiv f(y)\pmod Q$,
\begin{equation}
\begin{aligned}
\sum_{t\bmod S}e\!\left(\frac{\alpha f(t)+\beta t}{S}\right)
&=\sum_{y\bmod Q}e\!\left(\frac{\alpha f(y)+\beta y}{S}\right)
\sum_{v\bmod g}e\!\left(\frac{\beta v}{g}\right).
\end{aligned}
\end{equation}
The inner sum vanishes if $g\nmid\beta$. Otherwise write
$\alpha=g\alpha_0$, $\beta=g\beta_0$, and set
$F(t)=\alpha_0f(t)+\beta_0t$. As $\alpha\not\equiv0\pmod S$,
we have $Q>1$ and $\alpha_0\in\mathbb Z_p^\times$. Thus
\begin{equation}
\sum_{t\bmod S}e\!\left(\frac{\alpha f(t)+\beta t}{S}\right)
=g\sum_{t\bmod Q}e\!\left(\frac{F(t)}{Q}\right).
\end{equation}
Moreover,
\begin{equation}
\begin{aligned}
F'(t)&\equiv2\alpha_0a_2t+\alpha_0a_1+\beta_0\pmod p,\\
F''(t)&\equiv2\alpha_0a_2\not\equiv0\pmod p.
\end{aligned}
\end{equation}
Hensel's lemma gives a unique $\tau\in\mathbb Z_p$ with
$F'(\tau)=0$. Taylor expansion at $\tau$ then yields
\begin{equation}
F(\tau+v)=F(\tau)+\alpha_0a_2v^2\bigl(1+ph_\tau(v)\bigr),
\qquad h_\tau\in\mathbb Z_p[v].
\end{equation}
Since $p$ is odd, the $p$-adic square root
$\bigl(1+ph_\tau(v)\bigr)^{1/2}$ is defined for every
$v\in\mathbb Z_p$. The map
\begin{equation}
\psi(v)=v\bigl(1+ph_\tau(v)\bigr)^{1/2}
\end{equation}
is well defined modulo $Q$ and satisfies
$\psi(v)\equiv v\pmod p$ and $\psi'(v)\equiv1\pmod p$.
By Hensel's lemma, it is a bijection modulo $Q$. Consequently,
\begin{equation}
\begin{aligned}
\left|\sum_{t\bmod S}
e\!\left(\frac{\alpha f(t)+\beta t}{S}\right)\right|
&=g\left|\sum_{v\bmod Q}
e\!\left(\frac{\alpha_0a_2\psi(v)^2}{Q}\right)\right|\\
&=g\left|\sum_{v\bmod Q}
e\!\left(\frac{\alpha_0a_2v^2}{Q}\right)\right|\\
&=g\sqrt Q=\sqrt{Sg},
\end{aligned}
\end{equation}
by \eqref{frequency:quadratic-formula}.
\end{proof}

\section{Short second moments and local coordinates}\label{sec:short}

\subsection{The short moment formula}\label{short:reduction}

The following form of the reduction in \cite[Proposition~1]{MW} will be used in both counting arguments. We give the coefficient construction because its independence of the center character is needed later.
\begin{lemma}[Uniform short-moment reduction]\label{lem:short}
Let $p$ be a fixed odd prime, $q=p^n$, and let $p\le q_1<q$
be a power of $p$. Let
\begin{equation}\label{eq:short-parameters}
\begin{gathered}
D_h=q_1h,\qquad Q_h=\frac{q}{D_h},\qquad N_h=\frac{q}{D_h^2},\\
Q_0=Q_1=\frac{q}{q_1},\qquad N_0=N_1=\frac{q}{q_1^2},
\qquad  h=p^\ell\mid Q_0,\qquad \ell\ge0.
\end{gathered}
\end{equation}
For every $\varepsilon>0$, there are $O_\varepsilon(\log q)$
dyadic coefficient families $\{a_{N,h}(m)\}_h$, independent
of the primitive center $\chi\bmod q$, such that
\begin{equation}\label{eq:short-moment-reduction}
\begin{aligned}
S_2(\chi;q_1)\ll_{p,\varepsilon}q^\varepsilon q_1
\left(1+\max_{N}\sum_h Q_h^{-1/2}
\left|\sum_{\substack{0<|m|\le C_\varepsilon q^\varepsilon N_h\\p\nmid m}}
a_{N,h}(m)K_\chi(m;Q_h)\right|\right).
\end{aligned}
\end{equation}
Here $N$ ranges over the dyadic lengths in the approximate
functional equation, and
$$
|a_{N,h}(m)|\ll_\varepsilon q^\varepsilon,\qquad
h\ll_\varepsilon\frac{q^{1/2+\varepsilon}}{q_1},\qquad
Q_h\gg_\varepsilon q^{1/2-\varepsilon}
$$
whenever the corresponding coefficient system is nonzero.
For sufficiently large $q$, the kernels admit the decomposition
\eqref{eq:branches}. One may select a stationary sign and a fixed
unit class of $m\bmod p$ at a cost depending only on $p$.
\end{lemma}

\begin{proof}
By the approximate functional equation in Lemma \ref{AFE2}, we have
\begin{equation}\label{eq:short-afe}
\left|L\left(\frac12,\chi\psi\right)\right|^2
\ll (\log q)
\sum_{\substack{N\le q^{1/2+\varepsilon}\\N\ \mathrm{dyadic}}}
\left|\frac{1}{\sqrt{N}}\sum_n(\chi\psi)(n)V_N(n)\right|^2
+q^{-100},
\end{equation}
where $V_N$ is supported on $[N/2,2N]$,  $V_N^{(j)}(x) \ll_j N^{-j}$ for $j \in \mathbb{N}$.
By \eqref{intro:short} and \eqref{eq:short-afe}, we obtain
\begin{equation}\label{S_2 upper bound}
  S_2(\chi;q_1) \ll \frac{q^{\varepsilon}}{N}\sum_{\psi \bmod q_1}\left|\sum_n(\chi\psi)(n)V_N(n)\right|^2
+q^{-99},
\end{equation}
for some $N\le q^{1/2+\varepsilon}$.
By character orthogonality, we have 
$$
\begin{aligned}
&\frac1N\sum_{\psi \bmod q_1}
\left|\sum_n(\chi\psi)(n)V_N(n)\right|^2\\
&\qquad=\frac{\varphi(q_1)}N
\sum_{\substack{n_1\equiv n_2\pmod{q_1}\\p\nmid n_1n_2}}
\chi(n_1)\overline{\chi(n_2)}
V_N(n_1)\overline{V_N(n_2)}.
\end{aligned}
$$
The diagonal contribution is  $O(q_1)$. For $n_1=x+q_1b$,
$n_2=x$, $b>0$, let
$W_b(x)=V_N(x+q_1b)\overline{V_N(x)}$.
Pairing the positive and negative shifts, the off-diagonal
contribution is
\begin{equation}\label{eq:short-offdiagonal}
\frac{2\varphi(q_1)}N\operatorname{Re}
\sum_{1\le b\ll N/q_1}
\sum_{\substack{x\in\mathbb Z\\p\nmid x}}
W_b(x)\chi\left(1+\frac{q_1b}{x}\right).
\end{equation}
The character factor $\chi\left(1+\frac{q_1b}{x}\right)$ is periodic modulo $Q_0=q/q_1$.
The Poisson summation transforms \eqref{eq:short-offdiagonal} into
\begin{equation}\label{eq:short-poisson}
\frac{2\varphi(q_1)}{N\sqrt{Q_0}}\operatorname{Re}
\sum_{1\le b\ll N/q_1}\sum_{j\in\mathbb Z}
\widehat W_b(j/Q_0)K_\chi(j,b;Q_0),
\end{equation}
where $$K_\chi(j,b;Q_0)=\frac1{\sqrt{Q_0}}
\sum_{x\bmod Q_0}^{*}
e\left(\frac{(A/q_1)\log_p(1+q_1b/x)+jx}{Q_0}\right).$$
By \eqref{S_2 upper bound} through \eqref{eq:short-poisson}, we obtain
\begin{equation}\label{S_2 combine}
    S_2(\chi;q_1) \ll q^{\varepsilon}\left(q_1 + \frac{2\varphi(q_1)}{N\sqrt{Q_0}}\bigg|
\sum_{1\le b\ll N/q_1}\sum_{j\in\mathbb Z}
\widehat W_b(j/Q_0)K_\chi(j,b;Q_0) \bigg| \right).  
\end{equation}
For a nonempty $b$-range and sufficiently large $q$,
\begin{equation}\label{eq:short-exceptional}
1\le b\ll\frac N{q_1}<\frac{Q_0}{p},\qquad
v_p(b)\le v_p(Q_0)-2,\qquad K_\chi(0,b;Q_0)=0.
\end{equation}
By \cite[Lemma~6]{MW}, the exceptional unequal-valuation
cases require one of $j,b$ to be divisible by $Q_0$
and the other to have valuation $v_p(Q_0)-1$.
Thus every nonzero term satisfies $v_p(j)=v_p(b)$.
Writing
$$
j=hj',\qquad b=hb',\qquad h=p^{v_p(b)},\qquad
p\nmid j'b',
$$
and substituting $x=b'/u$, we obtain
\begin{equation}\label{eq:short-descent}
\begin{aligned}
K_\chi(hj',hb';Q_0)
&=\frac h{\sqrt{Q_0}}\sum_{u\bmod Q_h}^{*}
e\left(\frac{(A/D_h)\log_p(1+D_hu)+j'b'/u}{Q_h}\right)\\
&=h^{1/2}K_\chi(j'b';Q_h).
\end{aligned}
\end{equation}

Repeated integration by parts gives, for every fixed $J>0$,
$$
\frac{|\widehat W_b(u)|}{N}\ll_J(1+N|u|)^{-J}.
$$
Taking $J$ sufficiently large and using
$|K_\chi(j,b;Q_0)|\le\sqrt{Q_0}$, we may restrict
to $0<|j|\le q^\varepsilon Q_0/N$, with total error
$O(q^{-100})$. Hence, for $m=b'j'$,
\begin{equation}\label{eq:short-length}
\begin{gathered}
0<|m|\ll\frac{N}{q_1h}\frac{q^\varepsilon Q_0}{hN}
=q^\varepsilon N_h,\\
Q_h=\frac{Q_0}{h}\ge\frac{Q_0}{b}
\gg\frac qN\ge q^{1/2-\varepsilon},
\qquad h\le b\ll\frac{q^{1/2+\varepsilon}}{q_1}.
\end{gathered}
\end{equation}
For $p\nmid m$, define
\begin{equation}\label{eq:short-coefficients}
a_{N,h}(m)=
\sum_{\substack{b'>0,\ b'\mid m\\hb'\ll N/q_1}}
\frac{\widehat W_{hb'}(hm/(b'Q_0))}{N}
\mathbf 1_{\{|hm/b'|\le q^\varepsilon Q_0/N\}}.
\end{equation}
Then $a_{N,h}(m)$ is independent of $\chi$, and
$$
|a_{N,h}(m)|\ll d(|m|)\ll_\varepsilon q^\varepsilon.
$$
Grouping the terms according to $m=b'j'$ and
using \eqref{eq:short-coefficients}, we obtain
\begin{equation}\label{grouping the terms}
\begin{aligned}
&\sum_{1\le b\ll N/q_1}
\sum_{0<|j|\le q^\varepsilon Q_0/N}
\frac{\widehat W_b(j/Q_0)}{N}K_\chi(j,b;Q_0)\\
&\qquad=\sum_{h=p^\ell}h^{1/2}
\sum_{\substack{0<|m|\ll q^\varepsilon N_h\\p\nmid m}}
a_{N,h}(m)K_\chi(m;Q_h).
\end{aligned}
\end{equation}
Since $Q_0^{-1/2}h^{1/2}=Q_h^{-1/2}$,
\eqref{eq:short-moment-reduction} follows from
\eqref{S_2 combine} and \eqref{grouping the terms},
with the logarithmic factors absorbed into $q^\varepsilon$.
Finally, $Q_h\gg q^{1/2-\varepsilon}\ge p^2$ for
sufficiently large $q$. Splitting
$K_\chi=K_\chi^++K_\chi^-$ as in \eqref{eq:branches}
and partitioning $m\bmod p$ costs at most $2(p-1)$.
The finitely many remaining conductors are absorbed
into the implied constant.
\end{proof}

\subsection{Exact linearization of the stationary phase}
\label{sec:phase}

Lemma~\ref{lem:short} reduces the short second moment to
sums involving the complete kernels $K_A(m;Q)$. To study
their dependence on the character parameter within a fixed
coset, we introduce a local coordinate in which the
stationary phase becomes linear modulo $Q$. The following
proposition gives the precise form of this reduction.

\begin{proposition}[Local phase model]\label{prop:phase}
Let $p$ be an odd prime, $q=p^n$, and let
$$
D=p^d,\qquad R=p^s,\qquad Q=\frac qD,
\qquad q_2=DR<q,\qquad \Delta=\frac q{q_2}.
$$
Assume $d,s\ge1$ and $Q\ge p^2$, so that $Q=\Delta R$.
Fix an $H_{q_2}$-coset and a square-root branch
modulo $p$. Choose $c,B_0\in\mathbb Z_p^\times$ such that
$A\equiv cB_0^2\pmod\Delta$ throughout the coset.
On the fixed square-root branches, write
$A=cB^2$ and $m=cy^2$, with
$B\equiv B_0\pmod p$, where
$y=y(m)\in\mathbb Z_p^\times$ is independent of $A$.
Suppose that
\begin{equation}\label{eq:local-condition1}
5d\ge n+\mathbf 1_{p=5},
\end{equation}
\begin{equation}\label{eq:local-condition2}
2s\le n+d-\mathbf 1_{p=3},
\end{equation}
and
\begin{equation}\label{eq:local-condition3}
s\le 4d-2\mathbf 1_{p=3}.
\end{equation}
Then the parametrization
$B\equiv B_0+\Delta w\pmod Q$, with $w\bmod R$,
gives a bijection between the allowed parameters $A\bmod Q$
and $\mathbb Z/R\mathbb Z$.

Let
$$
\eta=-\frac{D^2}{24B_0^2}.
$$
Then $v_p(\eta)=2d-\mathbf 1_{p=3}\ge1$ and
$\eta^2\equiv0\pmod R$. There exist
$C_\pm(m)\in\mathbb Z_p^\times$, independent of $w$,
such that
\begin{equation}\label{eq:local-phase}
e\left(\frac{\Phi^\pm_{A,D}(m)}{Q}\right)
=
e\left(\frac{C_\pm(m)}{Q}\right)
e\left(\frac{\pm2cw(y+\eta y^3)}{R}\right).
\end{equation}
Moreover, on each fixed square-root branch,
the normalized stationary Gauss factors in
\eqref{eq:branches} are independent of $w$.
\end{proposition}

\begin{proof}
Write
$$
f_A(x)=\frac AD\log_p(1+Dx)+\frac mx.
$$
The stationary points of $f_A$ on the fixed
square-root branches are
$$
x_\pm=
\frac{Dm\pm\sqrt{D^2m^2+4Am}}{2A}.
$$
By \eqref{eq:branches}, the corresponding
stationary values are
$\Phi^\pm_{A,D}(m)=f_A(x_\pm)$.
For the subsequent expansion, let
$$
r^2=\frac mA,\qquad r=\frac yB,\qquad
z=\frac{Dr}{2}\in p\mathbb Z_p.
$$
Then
$$
x_\pm=r\left(z\pm\sqrt{1+z^2}\right),
$$
where $\sqrt{1+z^2}\equiv1\pmod p$.
Consequently, we have
\begin{equation*}
\begin{aligned}
1+Dx_\pm
  &=\left(\sqrt{1+z^2}\pm z\right)^2,\\
\log_p(1+Dx_\pm)
  &=\pm2\operatorname{arsinh}z,\\
\frac m{x_\pm}
  &=Ar\left(\pm\sqrt{1+z^2}-z\right),
\end{aligned}
\end{equation*}
where
$\operatorname{arsinh}z
=\log_p(z+\sqrt{1+z^2})$.
It follows that
$$
\Phi^\pm_{A,D}(m)
=
-\frac{Dm}{2}
\pm Ar\left(
\frac{\operatorname{arsinh}z}{z}
+\sqrt{1+z^2}
\right).
$$

Combining the binomial expansions of the two functions
in parentheses gives
$$
\frac{\operatorname{arsinh}z}{z}
+\sqrt{1+z^2}
=
2+\sum_{j\ge1}
\frac{2(-1)^{j-1}\binom{2j}{j}}
{4^j(4j^2-1)}z^{2j}.
$$
Since $z=Dr/2$, we obtain
\begin{equation}\label{eq:local-taylor}
\Phi^\pm_{A,D}(m)=
-\frac{Dm}{2}
\pm2Ar
\pm\sum_{j\ge1}c_jAD^{2j}r^{2j+1},
\end{equation}
with $$c_j
=
\frac{
2(-1)^{j-1}\binom{2j}{j}
}{16^j(4j^2-1)}.$$
In particular,
\begin{equation}\label{eq:local-coefficients}
c_1=\frac1{12},\qquad
c_2=-\frac1{320},\qquad
c_3=\frac1{3584}, \qquad
v_p(c_j)\ge-v_p(4j^2-1).
\end{equation}

We now show that the terms with $j\ge2$ can be
discarded modulo $Q$. Since $A,r\in\mathbb Z_p^\times$,
the second term in the tail satisfies
$$
v_p(c_2AD^4r^5)=4d-\mathbf 1_{p=5}.
$$
For $j\ge3$, the factors $2j-1$ and $2j+1$ are coprime
at every odd prime. Thus
\begin{equation*}
\begin{aligned}
v_p(4j^2-1)
&\le
\left\lfloor
\frac{\log(2j+1)}{\log 3}
\right\rfloor
\le 2j-4,\\
v_p(c_jAD^{2j}r^{2j+1})
&\ge 2jd-(2j-4)\\
&=4d+(2j-4)(d-1)
\ge4d.
\end{aligned}
\end{equation*}
Moreover, the lower bound
$2jd-\lfloor\log(2j+1)/\log3\rfloor$
tends to $+\infty$. Hence the series converges, and
\begin{equation*}
\sum_{j\ge2}c_jAD^{2j}r^{2j+1}
\in p^{4d-\mathbf 1_{p=5}}\mathbb Z_p.
\end{equation*}
By $Q=p^{n-d}$ and \eqref{eq:local-condition1}, we have
$4d-\mathbf 1_{p=5}\ge n-d=v_p(Q)$.
Therefore
\begin{equation}\label{eq:local-truncation}
\Phi^\pm_{A,D}(m)
\equiv
-\frac{Dm}{2}
\pm2Ar
\pm\frac{AD^2r^3}{12}
\pmod Q.
\end{equation}

We next introduce the local coordinate $w$.
Since $A=cB^2$ and $A\equiv cB_0^2\pmod\Delta$,
we have
$(B-B_0)(B+B_0)\equiv0\pmod\Delta$.
But $B+B_0\equiv2B_0\not\equiv0\pmod p$,
so $B+B_0$ is a unit and
$B\equiv B_0\pmod\Delta$. Thus
$B\equiv B_0+\Delta w\pmod Q$ for a unique
$w\bmod R$. Expanding $A=cB^2$ gives
\begin{equation}\label{eq:local-coordinate}
\frac{A-cB_0^2}{\Delta}
\equiv
2cB_0w+c\Delta w^2
\pmod R.
\end{equation}
To verify that this parametrization is bijective,
observe that for $w_1,w_2\bmod R$,
\begin{equation*}
\begin{aligned}
&\bigl(2cB_0w_1+c\Delta w_1^2\bigr)
-\bigl(2cB_0w_2+c\Delta w_2^2\bigr)\\
&\qquad=
c(w_1-w_2)
\bigl(2B_0+\Delta(w_1+w_2)\bigr).
\end{aligned}
\end{equation*}
The last factor is a $p$-adic unit. Hence the
right-hand side is divisible by $R$ if and only if
$w_1\equiv w_2\pmod R$. The map is therefore
injective, and hence bijective modulo $R$.

Using $r=y/B$ in \eqref{eq:local-truncation}, we obtain
$$
\Phi^\pm_{A,D}(m)
\equiv
-\frac{Dm}{2}
\pm2cBy
\pm\frac{cD^2y^3}{12B}
\pmod Q.
$$
We now use the exact identity
\begin{equation}\label{eq:local-reciprocal}
\frac1{B_0+\Delta w}
=
\frac1{B_0}
-\frac{\Delta w}{B_0^2}
+\frac{\Delta^2w^2}
{B_0^2(B_0+\Delta w)}.
\end{equation}
Substituting this into the preceding expression and
collecting the terms independent of $w$ and linear
in $w$ yields
\begin{equation}\label{eq:local-linearization}
\begin{aligned}
\Phi^\pm_{A,D}(m)
\equiv{}&
C_\pm(m)
\pm2c\Delta w
\left(y-\frac{D^2y^3}{24B_0^2}\right)\\
&\pm
\frac{cD^2\Delta^2w^2y^3}
{12B_0^2(B_0+\Delta w)}
\pmod Q,
\end{aligned}
\end{equation}
where
\begin{equation}\label{eq:local-constant}
C_\pm(m)
=
-\frac{Dm}{2}
\pm2cB_0y
\pm\frac{cD^2y^3}{12B_0}.
\end{equation}

The last term in \eqref{eq:local-linearization}
vanishes modulo $Q$. Indeed, since $B_0$ and
$B_0+\Delta w$ are units,
\begin{equation*}
\begin{aligned}
v_p\left(
\frac{cD^2\Delta^2w^2y^3}
{12B_0^2(B_0+\Delta w)}
\right)
&\ge
2d+2(n-d-s)-\mathbf 1_{p=3}\\
&=2n-2s-\mathbf 1_{p=3}\\
&\ge n-d=v_p(Q),
\end{aligned}
\end{equation*}
where the last inequality follows from
\eqref{eq:local-condition2}. Thus
\begin{equation*}
\Phi^\pm_{A,D}(m)
\equiv
C_\pm(m)
\pm2c\Delta w(y+\eta y^3)
\pmod Q,
\qquad
\eta=-\frac{D^2}{24B_0^2}.
\end{equation*}

Furthermore,
$v_p(\eta)=2d-\mathbf 1_{p=3}\ge1$, and
\eqref{eq:local-condition3} gives
$$
2v_p(\eta)
=
4d-2\mathbf 1_{p=3}
\ge s.
$$
Consequently, $\eta\in p\mathbb Z_p$ and
$\eta^2\equiv0\pmod R$.
Also, the first and third terms in
\eqref{eq:local-constant} are divisible by $p$.
Hence
$$
C_\pm(m)
\equiv\pm2cB_0y\not\equiv0\pmod p,
$$
so $C_\pm(m)\in\mathbb Z_p^\times$.

Finally, using $Q=\Delta R$, the preceding
congruence gives
\begin{equation*}
\begin{aligned}
e\left(\frac{\Phi^\pm_{A,D}(m)}{Q}\right)
&=
e\left(\frac{C_\pm(m)}{Q}\right)
e\left(\frac{\pm2c\Delta w(y+\eta y^3)}{Q}\right)\\
&=
e\left(\frac{C_\pm(m)}{Q}\right)
e\left(\frac{\pm2cw(y+\eta y^3)}{R}\right),
\end{aligned}
\end{equation*}
which proves \eqref{eq:local-phase}.

It remains to check the stationary Gauss factors.
Differentiating $f_A$ twice and using the stationary
equation $m=Ax_\pm^2/(1+Dx_\pm)$, we find
\begin{equation*}
\frac12 f_A''(x_\pm)
=
\frac{A(2+Dx_\pm)}
{2x_\pm(1+Dx_\pm)^2}
\equiv\frac A{x_\pm}\pmod p.
\end{equation*}
On the chosen branch,
$A\equiv cB_0^2\pmod p$ and
$x_\pm\equiv\pm y/B_0\pmod p$.
Thus the Hessian modulo $p$ is independent of $w$.
Since the normalized quadratic Gauss factor depends
only on this Hessian and the parity of $v_p(Q)$,
it is also independent of $w$. This completes the
proof of Proposition~\ref{prop:phase}.
\end{proof}

\begin{remark}
The condition $\eta^2\equiv0\pmod R$ is not needed
for the exact linearization of the stationary phase.
Indeed, \eqref{eq:local-condition1} and
\eqref{eq:local-condition2} already suffice to
establish \eqref{eq:local-phase}. The remaining
condition \eqref{eq:local-condition3} ensures
$\eta^2\equiv0\pmod R$, which is used in
Section~\ref{sec:finite} to obtain the exact inverse
$$
h^{-1}(x)\equiv x-\eta x^3\pmod R,
\qquad h(x)=x+\eta x^3,
$$
and to control congruences among the frequencies
$\xi_m=r(m)+\eta r(m)^3$.
Thus the nilpotence condition is imposed for the
subsequent congruence-spacing and higher-energy
arguments, rather than for the phase linearization
itself.
\end{remark}

\section{Large values of finite-ring square-root sums}\label{sec:finite}

The proof of Theorem~\ref{thm:finite} relies on two
estimates: a higher congruence-energy bound
for the large-value set and a coefficient-free
geometric bound obtained by Poisson summation.
The former follows from a spacing estimate
for the perturbed square-root frequencies,
while the latter exploits an averaged estimate
for the resulting complete sums with
quadratic--quartic phases. We establish these
ingredients separately and combine them
in the final subsection.

Fix $R=p^s$, $c\in\mathbb Z_p^\times$, and $a\in\mathbb F_p^\times$. We fix a lift of $a$ to $\mathbb Z_p^\times$,
still denoted by $a$. 
For $m\equiv ca^2\pmod p$, Hensel's lemma defines a unique root
$r(m)\in\mathbb Z_p^\times$ by $r(m)^2=m/c$ and $r(m)\equiv a\pmod p$.
We assume
\begin{equation}\label{finite:setup}
2\le N\le R\le N^{4/3},\qquad
\eta\in p\mathbb Z_p,\quad \eta^2\equiv0\pmod R,\qquad
\xi_m=r(m)+\eta r(m)^3.
\end{equation}
For an ordinary integer interval $I$ containing at most $N$ integers, omit
$m\not\equiv ca^2\pmod p$ and set
\begin{equation}\label{eq:finite-poly}
\begin{gathered}
\mathcal D(w)=\sum_{m\in I}a_m e\!\left(\frac{w\xi_m}{R}\right),
\qquad |a_m|\le1,\\
Z>0,\quad |\mathcal D(w)|\ge Z\quad (w\in\mathcal W),
\end{gathered}
\end{equation}
where $\mathcal W\subset\mathbb Z/R\mathbb Z$. 
Write 
\begin{equation}\label{eq:nilpotence-params}
E=p^{\min\{v_p(\eta),s\}},\qquad \Upsilon=R/E.
\end{equation}
The map $x\mapsto cx^2$ bijects the class $x\equiv a\pmod p$ with
the class $m\equiv ca^2\pmod p$ modulo $R$.  Moreover,
$$
h(x)=x+\eta x^3,\qquad
h^{-1}(z)\equiv z-\eta z^3\pmod R,
$$
since $\eta^2\equiv0\pmod R$. Consequently $m\mapsto\xi_m$
is injective modulo $R$ on the chosen branch.

\begin{theorem}[Finite-ring large values]\label{thm:finite}
Under \eqref{finite:setup}--\eqref{eq:nilpotence-params}, assume
$RE\ge N^2$.  For every fixed integer $k\ge2$ and $\varepsilon>0$,
\begin{equation}\label{eq:finite-main}
\begin{aligned}
|\mathcal W|\ll_{p,k,\varepsilon}R^\varepsilon\bigg(
&\frac{N^2}{Z^2}
+\frac{N^{3-2/k}R^{(k+1)/(2k)}\Upsilon^{1/(4k)}}{Z^{4-2/k}}\\
&+\frac{N^{k+1/2}R^{(k+1)/2}\Upsilon^{1/4}}{Z^{2k+1}}
\bigg).
\end{aligned}
\end{equation}
\end{theorem}

We first establish spacing and energy bounds, then average the Poisson
kernel over its dual frequency. The method is related to
\cite[Sections~4--6]{DFP}.

\subsection{Spacing of square root differences}\label{finite:spacing}

\begin{lemma}\label{lem:spacing}
Assume $2\le N\le R$ and the branch and
nilpotence conditions in \eqref{finite:setup}. For
$\ell\not\equiv0\pmod R$, we have
\begin{equation}\label{eq:spacing-fine}
\begin{aligned}
C(\ell)&:=\#\{m_1,m_2\in I:
\xi_{m_2}-\xi_{m_1}\equiv\ell\pmod R\}\\
&\ll_p\frac{N^2}{R}+
\sqrt{\frac{R}{(\ell,R)}}\log^2(2R).
\end{aligned}
\end{equation}
Furthermore, if we assume $R\le N^{4/3}$, then
\begin{equation}\label{eq:spacing-uniform}
C(\ell)\ll_{p,\varepsilon}R^\varepsilon N^2/R.
\end{equation}
\end{lemma}
\begin{proof}
If $p\nmid\ell$, then $C(\ell)=0$. Since $r(m_1)\equiv r(m_2) \equiv a \pmod p$ and $\eta \equiv 0 \pmod p$, we have 
$\xi_{m_1}\equiv \xi_{m_2} \equiv a \pmod p$. Hence we assume $p\mid\ell $.
Write $d=(\ell,R)$, $S=R/d$, and $\ell=d\ell_0$.
By $\xi_m\equiv r(m)+\eta r(m)^3\pmod R$,
$\xi_{m_2}-\xi_{m_1}\equiv\ell\pmod R$ and $\eta^2\equiv0\pmod R$, we obtain
$$
r(m_2)\equiv r(m_1)+\ell
-\eta((r(m_1)+\ell)^3-r(m_1)^3)\pmod R.
$$
Squaring and using $\eta^2\equiv0\pmod R$,
we obtain
\begin{equation}\label{eq:root-difference}
r(m_2)^2-r(m_1)^2
\equiv\ell(2r(m_1)+\ell-2\eta P_3(r(m_1)))\pmod R,
\end{equation}
where $P_3(X)=3X^3+6\ell X^2+4\ell^2X+\ell^3.$ 
Moreover,
$$
d\mid(m_2-m_1), \qquad t=(m_2-m_1)/d, \qquad 0<|t|\le N/d.
$$
As $m_i=c\,r(m_i)^2$, division by $d$ in
\eqref{eq:root-difference} gives
$$
\frac{t}{c\ell_0}\equiv2r(m_1)+\ell-2\eta P_3(r(m_1))\pmod S.
$$
Define $z(t)=\frac12\left(\frac{t}{c\ell_0}-\ell\right)$. Because $\eta^2=0\pmod S$,
$$
\begin{aligned}
r(m_1)-\eta P_3(r(m_1))&\equiv z(t)\pmod S\\
&\Longleftrightarrow\quad
r(m_1)\equiv z(t)+\eta P_3(z(t))\pmod S.
\end{aligned}
$$
The inverse is unique since $1-\eta P_3'(r(m_1))\equiv1\pmod p$.
Squaring this identity yields
\begin{equation}\label{spacing:solve}
m_1\equiv f(t):=c\bigl(z(t)^2+2\eta z(t)P_3(z(t))\bigr)\pmod S.
\end{equation}
Here
$f(t)\equiv\frac{t^2}{4c\ell_0^2} \pmod p ,$
and $(4c\ell_0^2,p)=1$.
Thus $f(t)$ satisfies Lemma~\ref{lem:gauss}. 
The map $(m_1,m_2)\mapsto(t,m_1)$ is
injective. Dropping the condition $m_2\in I$
and the residue restriction on $m_1$, we have
$$
C(\ell)\le
\#\{(t,m)\in\mathbb Z\times I:0<|t|\le N/d,
\ m\equiv f(t)\pmod S\}.
$$
If $d>N$, then $C(\ell)=0$. If $S\le N$, then
\begin{equation}\label{eq:spacing-large-d}
C(\ell)\ll\frac Nd\left(1+\frac NS\right)
\ll\frac{N^2}{R}.
\end{equation}
It remains to consider $S>N$. 
Split the positive and negative $t$-ranges
into ordinary integer intervals $J$, each
containing at most $N/d<S$ integers. 
With $\widehat{1_I}(u)=\sum_{m\in I}e(-um/S)$,
orthogonality yields
\begin{equation}\label{spacing:incidence-inversion}
\#\{(t,m)\in J\times I:m\equiv f(t)\pmod S\}
=\frac1S\sum_{u\bmod S}\widehat{1_I}(u)
\sum_{t\in J}e\!\left(\frac{uf(t)}S\right).
\end{equation}
The term $u=0$ is $|I||J|/S\ll N^2/R$.
For $u\not\equiv0\pmod S$, Fourier completion gives
$$
\sum_{t\in J}e\!\left(\frac{uf(t)}S\right)
=\frac1S\sum_{b\bmod S}\widehat{1_J}(b)
\sum_{t\bmod S}e\!\left(\frac{uf(t)+bt}{S}\right).
$$
Since $S^{-1}\sum_b|\widehat{1_J}(b)|\ll\log(2S)$,
Lemma~\ref{lem:gauss} gives
\begin{equation}\label{eq:incomplete-f}
\left|\sum_{t\in J}e\!\left(\frac{uf(t)}S\right)\right|
\ll_p\sqrt{S(u,S)}\log(2S).
\end{equation}
For $0<|u|\le S/2$,
\begin{equation}
|\widehat{1_I}(u)|\ll\min(N,S/|u|).
\end{equation}
Writing $S=p^\nu$ and grouping
$u=p^jb$, $p\nmid b$, gives
\begin{equation}\label{eq:weighted-fourier}
\sum_{\substack{u\bmod S\\u\ne0}}
|\widehat{1_I}(u)|\sqrt{(u,S)}
\ll S\log(2S)\sum_{0\le j<\nu}p^{-j/2}
\ll_p S\log(2S).
\end{equation}
Combining \eqref{eq:spacing-large-d}--\eqref{eq:weighted-fourier},
$$
C(\ell)\ll_p\frac{N^2}{R}+\sqrt S\log^2(2S)
\le\frac{N^2}{R}+\sqrt{\frac{R}{(\ell,R)}}\log^2(2R).
$$
Finally, we have $\sqrt R\le N^2/R$ if $R\le N^{4/3}$.
\end{proof}

\subsection{Higher congruence energies}\label{sec:energy}
For $D_0\mid R$, let
\begin{equation}\label{eq:energy-def}
\mathcal E_{2k}(\mathcal W;D_0;R)
=\#\left\{\boldsymbol w\in\mathcal W^{2k}:
\sum_{i=1}^kw_i-\sum_{i=k+1}^{2k}w_i
\equiv0\pmod{R/D_0}\right\}.
\end{equation}
Write $\mathcal E_{2k}(D_0)$ when $\mathcal W,R$ are fixed.
Each class modulo $R/D_0$ contains at most $D_0$ elements
of $\mathcal W$; thus $\mathcal E_2(D_0)\le D_0|\mathcal W|$.

\begin{lemma}[Energy estimate]\label{lem:energy}
For every fixed integer $k\ge2$, we have
\begin{equation}\label{eq:energy-bound}
\mathcal E_{2k}(D_0)\ll_{p,k,\varepsilon}R^\varepsilon
\left(\left(\frac NZ\right)^{4k-4}\mathcal E_2(D_0)
+\frac{ND_0}{Z^2}|\mathcal W|^{2k-1}\right).
\end{equation}
\end{lemma}
\begin{proof}
For each tuple counted by $\mathcal E_{2k}(D_0)$, insert
$Z^2\le|\mathcal D(w_{2k})|^2$, then enlarge the $w_{2k}$-sum
from $\mathcal W$ to its entire residue class modulo $R/D_0$.
For $\Delta(\boldsymbol w)=\sum_{i\le k}w_i-\sum_{k<i<2k}w_i$,
the resulting inner sum is
$$
\sum_{\substack{w\bmod R\\w\equiv \Delta(\boldsymbol w)\;(R/D_0)}}
e\!\left(\frac{w(\xi_{m_1}-\xi_{m_2})}{R}\right)
=D_0\,\mathbf1_{\{D_0\mid\xi_{m_1}-\xi_{m_2}\}}
e\!\left(\frac{\Delta(\boldsymbol w)(\xi_{m_1}-\xi_{m_2})}{R}\right).
$$

Let
$
F(\ell)=\sum_{w\in\mathcal W}
e\left(\frac{w\ell}{R}\right)
$
denote the finite Fourier transform of
$\mathbf 1_{\mathcal W}$. Moreover,
$$
\sum_{w_1,\ldots,w_{2k-1}\in\mathcal W}
e\!\left(\frac{\Delta(\boldsymbol w)(\xi_{m_1}-\xi_{m_2})}{R}\right)
=F(\xi_{m_1}-\xi_{m_2})^k
\overline{F(\xi_{m_1}-\xi_{m_2})}^{\,k-1}.
$$
Expanding $|\mathcal D(w)|^2$ (with the necessary conjugates), we have
\begin{equation}\label{energy:exact-positive}
\begin{aligned}
Z^2\mathcal E_{2k}(D_0)
&\le D_0\!\sum_{\substack{m_1,m_2\in I\\D_0\mid\xi_{m_1}-\xi_{m_2}}}
a_{m_1}\overline{a_{m_2}}
F(\xi_{m_1}-\xi_{m_2})^k
\overline{F(\xi_{m_1}-\xi_{m_2})}^{\,k-1}\\
&\le D_0\!\sum_{\substack{m_1,m_2\in I\\D_0\mid\xi_{m_1}-\xi_{m_2}}}
|F(\xi_{m_1}-\xi_{m_2})|^{2k-1}.
\end{aligned}
\end{equation}
The first sum is real and nonnegative by its original construction.
The injectivity of $m\mapsto\xi_m\pmod R$ and $N\le R$ imply
that $\xi_{m_1}-\xi_{m_2}\equiv0\pmod R$ precisely when $m_1=m_2$.
Its contribution is $\le ND_0|\mathcal W|^{2k-1}$.
For $\xi_{m_1}-\xi_{m_2}\not\equiv0\pmod R$,
Lemma~\ref{lem:spacing} gives
$C(\xi_{m_1}-\xi_{m_2})\ll R^\varepsilon N^2/R$.
Orthogonality gives, for $j=1,k$,
\begin{equation}\label{eq:energy-parseval}
\sum_{\substack{\ell\bmod R\\D_0\mid\ell}}|F(\ell)|^{2j}
=\frac R{D_0}\mathcal E_{2j}(D_0).
\end{equation}
H\"older with exponents $2k-2$ and $(2k-2)/(2k-3)$ gives
\begin{equation}\label{energy:holder-exponents}
\begin{aligned}
\sum_{\substack{\ell\bmod R\\D_0\mid\ell}}
|F(\ell)|^{2k-1}
&\le
\left(
\sum_{\substack{\ell\bmod R\\D_0\mid\ell}}
|F(\ell)|^2
\right)^{1/(2k-2)}
\left(
\sum_{\substack{\ell\bmod R\\D_0\mid\ell}}
|F(\ell)|^{2k}
\right)^{(2k-3)/(2k-2)}
\nonumber\\
&=
\frac{R}{D_0}
\mathcal E_2(D_0)^{1/(2k-2)}
\mathcal E_{2k}(D_0)^{(2k-3)/(2k-2)}.
\end{aligned}
\end{equation}
In particular,
$$
D_0\cdot\frac{N^2}{R}\cdot\frac R{D_0}=N^2.
$$
Therefore
\begin{equation}\label{eq:energy-solve}
\begin{aligned}
Z^2\mathcal E_{2k}(D_0)\ll R^\varepsilon\bigl(&ND_0|\mathcal W|^{2k-1}
+N^2\mathcal E_2(D_0)^{1/(2k-2)}\\
&\times \mathcal E_{2k}(D_0)^{(2k-3)/(2k-2)}\bigr).
\end{aligned}
\end{equation}
Writing $X=\mathcal E_{2k}(D_0)$ and $Y=\mathcal E_2(D_0)$,
this yields either
$$
X\ll R^\varepsilon ND_0Z^{-2}|\mathcal W|^{2k-1}
$$
or
$$
Z^2X^{1/(2k-2)}\ll R^\varepsilon N^2Y^{1/(2k-2)}.
$$
Hence, \eqref{eq:energy-bound} is proved.
\end{proof}

\subsection{The Poisson kernel and conductor descent}\label{finite:poisson}
\begin{lemma}[Evaluation of the Poisson kernel]\label{lem:poisson-kernel}
Assume \eqref{finite:setup}, and write
\begin{equation}\label{finite:kernel-definition}
C(\ell,z)=\sum_{\substack{x\bmod R\\x\equiv a\;(p)}}
e\!\left(\frac{z(x+\eta x^3)-c\ell x^2}{R}\right).
\end{equation}
If $R/p\mid\ell$, then $C(\ell,z)=0$ unless $R/p\mid z$,
and $|C(\ell,z)|\le R/p$. Otherwise write
$d=(\ell,R)$, $S=R/d\ge p^2$, $\ell=db$, $p\nmid b$.
Then $C(\ell,z)=0$ unless $d\mid z$.  For $z=dz'$,
\begin{equation}\label{eq:C-stationary}
C(db,dz')=d\sqrt S\,\omega_{S,c,b}
 e\!\left(\frac{\Psi_b(z')}{S}\right)
\mathbf1_{\{z'\equiv2cba\;(p)\}}, 
\end{equation}
with 
$$\Psi_b(z')
=\frac{z'^2}{4cb}
+\frac{\eta z'^4}{8c^3b^3}, 
\qquad \omega_{S,c,b}=S^{-1/2}\sum_{y\bmod S}
e\!\left(\frac{-cby^2}{S}\right),
\qquad |\omega_{S,c,b}|=1.$$
The factor $\omega_{S,c,b}$ is independent of $z'$.
\end{lemma}

\begin{proof}
If $R/p\mid\ell$, the term $e(-c\ell x^2/R)=e(-c\ell a^2/R)$ is constant
for $x\equiv a\pmod p$. Since $h(x)=x+\eta x^3$ is a bijection
on this class, we have
\begin{equation*}
\sum_{\substack{y\bmod R\\y\equiv a\pmod p}}
e\left(\frac{zy}{R}\right)
=
\sum_{t\bmod R/p}
e\left(\frac{z(a+pt)}{R}\right)
=
e\left(\frac{za}{R}\right)
\sum_{t\bmod R/p}
e\left(\frac{zt}{R/p}\right).
\end{equation*}
By orthogonality of additive characters, we have 
$$
C(\ell,z)=e\!\left(\frac{-c\ell a^2}{R}\right)
\sum_{\substack{y\bmod R\\y\equiv a\pmod p}}e\!\left(\frac{zy}{R}\right)
=e\!\left(\frac{-c\ell a^2+za}{R}\right)
\sum_{t\bmod R/p}e\!\left(\frac{zt}{R/p}\right).
$$

Suppose $R/p\nmid\ell$. Write $d=(\ell, R)$, $\ell=db$, $R=dS$.
If $d\nmid z$, write
$h=v_p(z)<v_p(\ell)\le s-2 $, and $ z=p^h z_0$,
$ p\nmid z_0$,
and write $\delta=p^{s-h-1}$. For
$$
H(x)=z(x+\eta x^3)-c\ell x^2,
$$
a direct expansion, using $\eta\in p\mathbb Z_p$
and $v_p(\ell)>h$, gives
\begin{equation}
H(x+\delta)-H(x)
\equiv z\delta
=p^{s-1}z_0\pmod R.
\end{equation}
Since $\delta\equiv0\pmod p$, the translation
$x\mapsto x+\delta$ permutes the residue class
$x\equiv a\pmod p$. Consequently,
\begin{equation}
C(\ell,z)
=e\left(\frac{z_0}{p}\right)C(\ell,z).
\end{equation}
As $p\nmid z_0$, $e\left(\frac{z_0}{p}\right) \neq 1$.
Hence $C(\ell,z)=0$.

For $z=dz'$, conductor descent gives
$$
C(db,dz')
=d\sum_{\substack{x\bmod S\\x\equiv a\pmod p}}
e\!\left(\frac{G(x)}{S}\right),$$
where
$G(x)=z'(x+\eta x^3)-cbx^2.$
Let $t=z'/(2cb)$. Using 
$x^3-t^3=(x-t)(x^2+tx+t^2)$, we make the change of variables
$$
v=x-t-\eta t(x^2+tx+t^2).
$$
Since $\eta^2\equiv0\pmod S$ and $z'=2cbt$,
direct expansion yields
\begin{equation}
G(x)
\equiv \Psi_b(z')-cbv^2\pmod S,
\end{equation}
where  $$\Psi_b(z')
=\frac{z'^2}{4cb}
+\frac{\eta z'^4}{8c^3b^3}.$$
Moreover,
$$
v(x)-v(y)
=(x-y)\bigl(1-\eta t(x+y+t)\bigr).
$$
Thus $x\mapsto v$ is bijective modulo $S$, and
$v\equiv x-t\pmod p$. Consequently,
\begin{equation}
\begin{aligned}
\sum_{\substack{x\bmod S\\x\equiv a\;(p)}}
e\!\left(\frac{G(x)}{S}\right)
&\quad=
e\!\left(\frac{\Psi_b(z')}{S}\right)
\sum_{\substack{v\bmod S\\ v\equiv a-t\pmod p}}
e\!\left(\frac{-cbv^2}{S}\right)\\
&\quad=
\mathbf 1_{\{z'\equiv2cba\;(p)\}}
e\!\left(\frac{\Psi_b(z')}{S}\right)
\sum_{v\bmod S}e\!\left(\frac{-cbv^2}{S}\right).
\end{aligned}
\end{equation}
Indeed, since $S=p^\nu$ with $\nu\ge2$,
each nonzero residue class $v\equiv v_0\pmod p$,
where $v_0\in\mathbb F_p^\times$, contributes zero
under the translation $v\mapsto v+p^{\nu-1}$,
which introduces the factor
$e(-2cbv_0/p)\ne1$.
The quadratic Gauss formula now gives
\eqref{eq:C-stationary}.
\end{proof}

\subsection{Averaging the dual frequency}
\label{finite:frequency}

For $S=p^\nu\mid R$ and $p\nmid b$, write
\begin{equation}\label{eq:dual-phase}
\Psi_b(z)=\frac{z^2}{4cb}
+\eta\frac{z^4}{8c^3b^3},
\qquad
\Delta_j=
\sum_{i=1}^ku_i^j-\sum_{i=k+1}^{2k}u_i^j.
\end{equation}
All quantities $\Delta_j$ are understood modulo $S$.
Whenever a $p$-adic valuation of a residue class
modulo $S$ is used, it is truncated at $\nu$.
Define
\begin{equation}\label{eq:I-def}
I_b(\boldsymbol u)=
\sum_{x\bmod S}
e\left(
\frac{
\sum_{i=1}^k\Psi_b(x-u_i)
-\sum_{i=k+1}^{2k}\Psi_b(x-u_i)
}{S}
\right).
\end{equation}

\begin{lemma}[Complete-frequency average]
\label{lem:frequency}
Assume that $c\in\mathbb Z_p^\times$ and
$\eta^2\equiv0\pmod S$. Let
$$
D_0=\frac{S}{(\Delta_1,S)},
\qquad
e_0=\min\{v_p(\eta),\nu\}.
$$
If $D_0>S/p^{e_0}$, then
$I_b(\boldsymbol u)=0$ for every unit $b$.
Moreover, for $K\ge1$, we have
\begin{equation}\label{eq:I-average}
\sum_{\substack{1\le|b|\le K\\p\nmid b}}
|I_b(\boldsymbol u)|
\ll_p S\left(1+KD_0^{-1/2}\right).
\end{equation}
\end{lemma}

\begin{proof}
Expanding the quadratic and quartic terms in
\eqref{eq:I-def}, we find that, up to a term
independent of $x$, the phase is
\begin{equation}\label{eq:I-phase}
-\frac{c^2\Delta_1b^2+\eta\Delta_3}
{2c^3b^3}x
+\frac{3\eta\Delta_2}{4c^3b^3}x^2
-\frac{\eta\Delta_1}{2c^3b^3}x^3.
\end{equation}
The omitted constant contributes only a factor
of modulus one to $I_b(\boldsymbol u)$.

If $\Delta_1\equiv0\pmod S$, then $D_0=1$,
and the trivial bound $|I_b(\boldsymbol u)|\le S$
immediately gives \eqref{eq:I-average}.
Henceforth assume that
$\Delta_1\not\equiv0\pmod S$, and let
$$
j=v_p(\Delta_1)<\nu,
\qquad D_0=p^{\nu-j}.
$$

Suppose first that $j<e_0$.
Since $c$ and $b$ are units, we have
$v_p(c^2\Delta_1b^2)=j$, $v_p(\eta\Delta_3)\ge e_0>j.$ 
Thus the linear coefficient in \eqref{eq:I-phase}
has valuation exactly $j$. The quadratic
coefficient has valuation at least $e_0\ge j+1$,
while the cubic coefficient has valuation at least
$v_p(\eta)+j\ge j+1$.
Consequently, after removing the common factor
$p^j$ and descending to modulus $p^{\nu-j}$,
the reduced phase has the form
$$
\alpha x+ph(x),
\qquad
\alpha\in\mathbb Z_p^\times,\quad
h\in\mathbb Z_p[x].
$$
Lemma~\ref{lem:nonstationary-cancellation}
therefore gives
$I_b(\boldsymbol u)=0$.
Since $j<e_0$ is equivalent to
$D_0>S/p^{e_0}$, this proves the support
assertion.

We now suppose that $j\ge e_0$.
The condition $\eta^2\equiv0\pmod S$ implies
$2v_p(\eta)\ge\nu$, and hence
$
v_p(\eta\Delta_1)
\ge v_p(\eta)+e_0\ge\nu.
$
Thus the cubic term in \eqref{eq:I-phase}
vanishes modulo $S$.
Let
$
b_0=\min\{v_p(3\eta\Delta_2),\nu\}.
$
The remaining phase is quadratic. By
\eqref{frequency:quadratic-formula}, we have
\begin{equation*}
|I_b(\boldsymbol u)|
=
\sqrt{Sp^{b_0}}\,
\mathbf 1_{\{
p^{b_0}\mid c^2\Delta_1b^2+\eta\Delta_3
\}}.
\end{equation*}

If $b_0\le j$, then
$$
|I_b(\boldsymbol u)|
\le\sqrt{Sp^j}
=S D_0^{-1/2}.
$$
Summing over the $O(K)$ admissible frequencies
proves \eqref{eq:I-average} in this case.

It remains to consider $b_0>j$.
A nonzero Gauss sum requires
\begin{equation}\label{eq:exceptional-frequency}
c^2\Delta_1b^2+\eta\Delta_3
\equiv0\pmod{p^{b_0}}.
\end{equation}
Since
$v_p(c^2\Delta_1b^2)=j<b_0$, this forces
$v_p(\eta\Delta_3)=j$.
After division by $p^j$, the congruence
\eqref{eq:exceptional-frequency} becomes
$$
c^2\frac{\Delta_1}{p^j}b^2
+\frac{\eta\Delta_3}{p^j}
\equiv0\pmod{p^{b_0-j}}.
$$
This is a quadratic congruence in the unit $b$.
Its derivative is
$
2c^2\frac{\Delta_1}{p^j}b\in\mathbb Z_p^\times.
$
Hence every root modulo $p$ is simple and
lifts uniquely by Hensel's lemma. Since the
reduced equation has degree two, it has at
most two roots modulo $p^{b_0-j}$.

Consequently, the number of integers
$1\le|b|\le K$ satisfying
\eqref{eq:exceptional-frequency} is
$O(1+Kp^{-(b_0-j)})$. Therefore
\begin{equation}\label{eq:I-average-proof}
\begin{aligned}
\sum_{\substack{1\le|b|\le K\\p\nmid b}}
|I_b(\boldsymbol u)|
&\ll_p
\left(1+\frac K{p^{b_0-j}}\right)
\sqrt{Sp^{b_0}}\\
&\le
S+K\sqrt{Sp^j}
=
S\left(1+\frac K{\sqrt{D_0}}\right).
\end{aligned}
\end{equation}
Here we used $b_0\le\nu$ and
$D_0=p^{\nu-j}$.
The factor $3$ in the definition of $b_0$
accounts for the additional $p$-adic valuation
when $p=3$, so no separate restriction on $p$
is needed. This completes the proof.
\end{proof}

\subsection{Geometric estimate and proof of Theorem~\ref{thm:finite}}

\begin{lemma}[Geometric estimate]\label{lem:geometric}
Let $\phi\ge0$ be a fixed smooth function of compact
support, and let $M_0\in\mathbb R$. For arbitrary
coefficients $|\gamma_w|\le1$, define
\begin{equation}\label{geometric:kernel-norm}
\mathscr Q_N(\mathcal W;\gamma)
=
N\sum_{\substack{m\in\mathbb Z\\m\equiv ca^2\;(p)}}
\phi\left(\frac{m-M_0}{N}\right)
\left|\sum_{w\in\mathcal W}
\gamma_w e\left(\frac{w\xi_m}{R}\right)\right|^2.
\end{equation}
Under the phase and modulus hypotheses of
Theorem~\ref{thm:finite}, including $RE\ge N^2$,
for every fixed integer $k\ge2$,
\begin{equation}\label{eq:geometric}
\begin{aligned}
\mathscr Q_N(\mathcal W;\gamma)
\ll_{p,k,\phi}\;&
N^2|\mathcal W|\\
&+
NR^{(k+1)/(2k)}
|\mathcal W|^{1-1/(2k)}
\left(
\sum_{\substack{D_0\mid R\\D_0\le\Upsilon}}
D_0^{-1/2}\mathcal E_{2k}(D_0)
\right)^{1/(2k)}.
\end{aligned}
\end{equation}
The implied constant is independent of $M_0$, and
no large-value condition is needed here.
\end{lemma}

\begin{proof}
By Hensel's lemma, the map $x\mapsto cx^2$
is a bijection modulo $R$ from the class
$x\equiv a\pmod p$ to the class
$m\equiv ca^2\pmod p$.
Expanding the square in
\eqref{geometric:kernel-norm} and applying Poisson
summation in each residue class $m\bmod R$, we obtain
\begin{equation}\label{finite:poisson-exact}
\begin{aligned}
\mathscr Q_N(\mathcal W;\gamma)
=\frac{N^2}{R}\sum_{\ell\in\mathbb Z}
&\widehat\phi\left(-\frac{\ell N}{R}\right)
e\left(\frac{\ell M_0}{R}\right)\\
&\times
\sum_{w,w'\in\mathcal W}
\gamma_w\overline{\gamma_{w'}}
C(\ell,w-w').
\end{aligned}
\end{equation}
The factor involving $M_0$ has modulus one, so it
does not affect any of the subsequent upper bounds.

We first consider the terminal frequencies
$R/p\mid\ell$. By Lemma~\ref{lem:poisson-kernel},
$C(\ell,w-w')$ vanishes unless
$R/p\mid w-w'$, and
$|C(\ell,w-w')|\le R/p$.
For each $w\in\mathcal W$, there are at most
$p$ admissible choices of $w'\in\mathcal W$.
Writing $\ell=jR/p$ and using the rapid decay of
$\widehat\phi$, their total contribution is
\begin{equation}\label{finite:terminal}
\begin{aligned}
\frac{N^2}{R}\frac Rp\,(p|\mathcal W|)
\sum_{j\in\mathbb Z}
\left(1+\frac{|j|N}{p}\right)^{-3}
\ll_{p,\phi}N^2|\mathcal W|.
\end{aligned}
\end{equation}
Here we used $N\ge1$.

We now turn to the nonterminal frequencies.
Let
$$
d=(\ell,R),\qquad
S=R/d\ge p^2,\qquad
\ell=db,\quad p\nmid b.
$$
By \eqref{eq:C-stationary}, the kernel vanishes unless
$w\equiv w'\pmod d$. Fix such a residue class and write
$$
w=w_0+du,\qquad
w'=w_0+dv,\qquad u,v\bmod S.
$$
The remaining support condition is
$
u-v\equiv2cba\pmod p.
$
Thus, after fixing $b\bmod p$ and $u\bmod p$,
the residue class $v\bmod p$ is uniquely determined.

Denote the corresponding scaled subsets by
$U,V\subseteq\mathbb Z/S\mathbb Z$, and let
$
\alpha_u=\gamma_{w_0+du}
$,
$
\beta_v=\overline{\gamma_{w_0+dv}}.
$
Then $|\alpha_u|,|\beta_v|\le1$.
Define
$$
H_b(x)=
\sum_{v\in V}\beta_v
e\left(\frac{\Psi_b(x-v)}S\right).
$$
By \eqref{eq:C-stationary}, the Gauss factor is
independent of $u-v$ and has modulus one.
Consequently, the contribution from these fixed
residue classes is bounded by
$$
\frac{N^2}{\sqrt S}
\sum_{\substack{b\in\mathbb Z\\
p\nmid b\\b\equiv b_*\;(p)}}
\left|\widehat\phi\left(-\frac{Nb}{S}\right)\right|
\sum_{u\in U}|H_b(u)|,
$$
where $b_*\not\equiv0\pmod p$ is fixed.
Denote this upper bound by $\mathcal T_d(U,V)$.

Fix $L=3$ and let
$$
\varpi(n)=\left(1+\frac{Nn}{S}\right)^{-3},
\qquad
\mathcal B=\sum_{n\ge1}\varpi(n).
$$
The rapid decay of $\widehat\phi$ gives
$$
\left|\widehat\phi\left(-\frac{Nb}{S}\right)\right|
\ll_\phi\varpi(|b|).
$$
Since
$\sum_{b\in\mathbb Z,\ p\nmid b}\varpi(|b|)
\le2\mathcal B$,
weighted H\"older's inequality with exponent $2k$,
followed by extending the nonnegative sum over
$u\in U$ to all $x\bmod S$, yields
\begin{equation}\label{finite:holder-full}
\mathcal T_d(U,V)
\ll_{\phi}\;
\frac{N^2}{\sqrt S}
(\mathcal B|U|)^{1-1/(2k)}
\times
\left(
\sum_{\substack{b\in\mathbb Z\\p\nmid b}}
\varpi(|b|)
\sum_{x\bmod S}|H_b(x)|^{2k}
\right)^{1/(2k)}.
\end{equation}
Here and below, dropping restrictions on the unit
frequencies only enlarges the nonnegative sums.

We next apply Lemma~\ref{lem:frequency} in a weighted
form. For $\boldsymbol v\in V^{2k}$, let
$$
\Delta_1(\boldsymbol v)
=
\sum_{i=1}^kv_i-\sum_{i=k+1}^{2k}v_i,
\qquad
D_*(\boldsymbol v)
=
\frac{S}{(\Delta_1(\boldsymbol v),S)}.
$$
If $D_*(\boldsymbol v)>S/p^{e_0}$, where
$e_0=\min\{v_p(\eta),v_p(S)\}$, then
$I_b(\boldsymbol v)=0$ for every $p\nmid b$.
Otherwise, summation by parts and
\eqref{eq:I-average} give
\begin{equation}\label{finite:abel-weighted}
\begin{aligned}
\sum_{\substack{b\in\mathbb Z\\p\nmid b}}
\varpi(|b|)|I_b(\boldsymbol v)|
&\quad=
\sum_{n\ge1}
\bigl(\varpi(n)-\varpi(n+1)\bigr)
\sum_{\substack{1\le|b|\le n\\p\nmid b}}
|I_b(\boldsymbol v)|\\
&\quad\ll_p
S\sum_{n\ge1}
\bigl(\varpi(n)-\varpi(n+1)\bigr)
\left(1+\frac n{\sqrt{D_*(\boldsymbol v)}}\right)\\
&\quad\ll_p
S\left(
\varpi(1)+\frac{\mathcal B}{\sqrt{D_*(\boldsymbol v)}}
\right).
\end{aligned}
\end{equation}
In the last step we used
$$
\sum_{n\ge1}(\varpi(n)-\varpi(n+1))
=\varpi(1), \qquad
\sum_{n\ge1}n(\varpi(n)-\varpi(n+1))
=\mathcal B.
$$

Expanding the even power gives
$$
\sum_{x\bmod S}|H_b(x)|^{2k}
=
\sum_{\boldsymbol v\in V^{2k}}
\left(\prod_{i=1}^k\beta_{v_i}\right)
\left(\prod_{i=k+1}^{2k}\overline{\beta_{v_i}}\right) I_b(\boldsymbol v).
$$
Taking absolute values and using
$|\beta_v|\le1$, we may sum the weighted estimate
\eqref{finite:abel-weighted} over all tuples.
Partition the tuples according to their exact divisor
$D_*(\boldsymbol v)$. If $D_*(\boldsymbol v)=D_0$,
then
$S/D_0\mid\Delta_1(\boldsymbol v)$.
Hence their number is at most
$\mathcal E_{2k}(V;D_0;S)$.
We conclude that
\begin{equation}\label{finite:tuple-bound}
\begin{aligned}
&\sum_{\substack{b\in\mathbb Z\\p\nmid b}}
\varpi(|b|)
\sum_{x\bmod S}|H_b(x)|^{2k}\\
&\qquad\ll_p
S\sum_{\substack{D_0\mid S\\D_0\le S/p^{e_0}}}
\left(\varpi(1)+\mathcal B D_0^{-1/2}\right)
\mathcal E_{2k}(V;D_0;S).
\end{aligned}
\end{equation}

Since $\varpi$ is positive and decreasing, integral
comparison gives
\begin{equation}\label{eq:L-sum}
\frac{S}{2N}\left(1+\frac NS\right)^{-2}
\le\mathcal B\le\frac{S}{2N},\qquad
\frac{\varpi(1)}{\mathcal B}
\le\frac{2N}{S+N}
\le\frac{2N}{S}.
\end{equation}

Recall that $S=p^\nu$, $R=p^s$ with $\nu\le s$,
and write $E=(R,\eta)$, so that $\Upsilon=R/E$.
Then
$$
\frac{S}{p^{e_0}}
=
p^{\max\{\nu-v_p(\eta),0\}}
\le
p^{\max\{s-v_p(\eta),0\}}
=\frac RE=\Upsilon.
$$
Consequently, for every $D_0$ occurring in
\eqref{finite:tuple-bound},
$$
\frac{\varpi(1)\sqrt{D_0}}{\mathcal B}
\ll
\frac NS\sqrt{\frac RE}
=
\frac{dN}{\sqrt{RE}}
\le d,
$$
where the last inequality follows from $RE\ge N^2$.
It follows that
$$
\varpi(1)+\mathcal B D_0^{-1/2}
\ll
\mathcal B(1+d)D_0^{-1/2}.
$$
Combining this with
\eqref{finite:holder-full},
\eqref{finite:tuple-bound}, and
\eqref{eq:L-sum}, we obtain
\begin{equation}\label{eq:dyadic-L}
\begin{aligned}
\mathcal T_d(U,V)
\ll_{p,k,\phi}\;&
NS^{(k+1)/(2k)}
(1+d)^{1/(2k)}
|U|^{1-1/(2k)}\\
&\times
\left(
\sum_{\substack{D_0\mid S\\D_0\le S/p^{e_0}}}
D_0^{-1/2}
\mathcal E_{2k}(V;D_0;S)
\right)^{1/(2k)}.
\end{aligned}
\end{equation}

It remains to combine the estimates over the
fixed residue-class configurations. For a tuple
in one scaled class, write $w_i=w_0+dv_i$.
The balanced signs give
\begin{equation}\label{finite:energy-embedding}
\begin{gathered}
\sum_{i=1}^kw_i-\sum_{i=k+1}^{2k}w_i
=
d\left(
\sum_{i=1}^kv_i-\sum_{i=k+1}^{2k}v_i
\right),\\
\frac S{D_0}\mid\Delta_1(\boldsymbol v)
\quad\Longrightarrow\quad
\frac R{D_0}\mid
\left(
\sum_{i=1}^kw_i-\sum_{i=k+1}^{2k}w_i
\right).
\end{gathered}
\end{equation}

For clarity, fix $b_*\bmod p$. A residue-class
configuration is determined by $w_0\bmod d$ and
$u_0\bmod p$. The corresponding sets are
$$
\begin{aligned}
U_\alpha
&=\{u\bmod S:
w_0+du\in\mathcal W,\ u\equiv u_0\pmod p\},\\
v_0&\equiv u_0-2cb_*a\pmod p,\\
V_\alpha
&=\{v\bmod S:
w_0+dv\in\mathcal W,\ v\equiv v_0\pmod p\}.
\end{aligned}
$$
For fixed $b_*$, the sets $U_\alpha$ are pairwise
disjoint under $u\mapsto w_0+du$. The same holds
for the sets $V_\alpha$, since the support condition
determines $v_0$ uniquely from $u_0$.
There are only $p-1$ possible unit classes $b_*$.
Thus the total multiplicity of these configurations
is $O_p(1)$.

Indexing the configurations by $\alpha$, let
$$
X_\alpha=|U_\alpha|,
\qquad
Y_\alpha=
\sum_{\substack{D_0\mid S\\D_0\le S/p^{e_0}}}
D_0^{-1/2}
\mathcal E_{2k}(V_\alpha;D_0;S).
$$
The disjointness and \eqref{finite:energy-embedding}
imply
$$
\sum_\alpha X_\alpha
\ll_p|\mathcal W|,\qquad
\sum_\alpha Y_\alpha
\ll_p
\sum_{\substack{D_0\mid R\\D_0\le\Upsilon}}
D_0^{-1/2}
\mathcal E_{2k}(\mathcal W;D_0;R).
$$
H\"older's inequality now yields
\begin{equation}\label{geometric:merge-holder}
\sum_\alpha
X_\alpha^{1-1/(2k)}Y_\alpha^{1/(2k)}
\le
\left(\sum_\alpha X_\alpha\right)^{1-1/(2k)}
\left(\sum_\alpha Y_\alpha\right)^{1/(2k)}.
\end{equation}

Finally, since $S=R/d$, we have
$
S^{(k+1)/(2k)}(1+d)^{1/(2k)}
\ll_k R^{(k+1)/(2k)}d^{-1/2}.
$
As $d$ ranges over divisors of the prime power $R$, we have
$
\sum_{d\mid R}d^{-1/2}
\le\sum_{j\ge0}p^{-j/2}\ll_p1.
$
Summing \eqref{eq:dyadic-L} over $d$, applying
\eqref{geometric:merge-holder}, and adding the
terminal contribution \eqref{finite:terminal}
proves \eqref{eq:geometric}.
\end{proof}

\begin{proof}[Proof of Theorem~\ref{thm:finite}]
If $|\mathcal W|=0$, the conclusion is immediate.
Choose a fixed nonnegative
$\phi\in C_c^\infty(\mathbb R)$ with
$\phi(t)\ge1$ for $|t|\le1$, and let $M_0$
be the midpoint of $I$. Since $I$ contains at most
$N$ consecutive integers,
$\phi((m-M_0)/N)\ge1$ for every $m\in I$.

Since $|\mathcal D(w)|\ge Z>0$ on $\mathcal W$,
let
$$
\gamma_w=
\frac{\overline{\mathcal D(w)}}{|\mathcal D(w)|},
\qquad |\gamma_w|=1.
$$
Then $\gamma_w\mathcal D(w)=|\mathcal D(w)|$.
By Cauchy--Schwarz and the definition of
$\mathscr Q_N$, we obtain
\begin{equation}\label{eq:duality}
\begin{aligned}
|\mathcal W|^2Z^2
&\le
\left|
\sum_{m\in I}a_m
\sum_{w\in\mathcal W}
\gamma_w e\left(\frac{w\xi_m}{R}\right)
\right|^2\\
&\le
\frac{\sum_{m\in I}|a_m|^2}{N}
\mathscr Q_N(\mathcal W;\gamma)
\le
\mathscr Q_N(\mathcal W;\gamma),
\end{aligned}
\end{equation}
where the last inequality follows from
$|a_m|\le1$ and $|I|\le N$.

By Lemma~\ref{lem:energy} and
$\mathcal E_2(D_0)\le D_0|\mathcal W|$, we have
$$
\begin{aligned}
D_0^{-1/2}\mathcal E_{2k}(D_0)
\ll_{p,k,\varepsilon}R^\varepsilon D_0^{1/2}
\left(
\left(\frac NZ\right)^{4k-4}|\mathcal W|
+\frac N{Z^2}|\mathcal W|^{2k-1}
\right).
\end{aligned}
$$
Since $D_0$ ranges over powers of $p$,
$
\sum_{\substack{D_0\mid R\\D_0\le\Upsilon}}
D_0^{1/2}\ll_p\Upsilon^{1/2}.
$
It follows that
\begin{equation}\label{eq:energy-insert}
\sum_{\substack{D_0\mid R\\D_0\le\Upsilon}}
D_0^{-1/2}\mathcal E_{2k}(D_0)
\ll_{p,k,\varepsilon}
R^\varepsilon\Upsilon^{1/2}
\left(
\left(\frac NZ\right)^{4k-4}|\mathcal W|
+\frac N{Z^2}|\mathcal W|^{2k-1}
\right).
\end{equation}

Combining \eqref{eq:duality},
Lemma~\ref{lem:geometric}, and
\eqref{eq:energy-insert}, and using
$(A+B)^{1/(2k)}
\le A^{1/(2k)}+B^{1/(2k)}$ for $A,B\ge0$,
we obtain
\begin{equation}\label{finite:three-terms-unsolved}
\begin{aligned}
|\mathcal W|^2Z^2\ll_{p,k,\varepsilon}R^\varepsilon
\bigg(&N^2|\mathcal W|\\
&+
N^{3-2/k}R^{(k+1)/(2k)}
\Upsilon^{1/(4k)}
Z^{-2+2/k}|\mathcal W|\\
&+
N^{1+1/(2k)}R^{(k+1)/(2k)}
\Upsilon^{1/(4k)}
Z^{-1/k}|\mathcal W|^{2-1/k}
\bigg).
\end{aligned}
\end{equation}

Dividing by $|\mathcal W|Z^2$, we find
$$
\begin{aligned}
|\mathcal W|\ll_{p,k,\varepsilon}R^\varepsilon\bigg(
&\frac{N^2}{Z^2}
+
\frac{
N^{3-2/k}R^{(k+1)/(2k)}
\Upsilon^{1/(4k)}
}{Z^{4-2/k}}\\
&+
\frac{
N^{1+1/(2k)}R^{(k+1)/(2k)}
\Upsilon^{1/(4k)}
}{Z^{2+1/k}}
|\mathcal W|^{1-1/k}
\bigg).
\end{aligned}
$$
We use the elementary inequality
$$
x\le A+Bx^{1-1/k}
\quad\Longrightarrow\quad
x\ll_k A+B^k,
\qquad A,B\ge0.
$$
Applying it to the preceding estimate, 
we obtain
$$
\begin{aligned}
|\mathcal W|\ll_{p,k,\varepsilon}R^\varepsilon\bigg(
&\frac{N^2}{Z^2}
+
\frac{
N^{3-2/k}R^{(k+1)/(2k)}
\Upsilon^{1/(4k)}
}{Z^{4-2/k}}\\
&+
\frac{
N^{k+1/2}R^{(k+1)/2}
\Upsilon^{1/4}
}{Z^{2k+1}}
\bigg).
\end{aligned}
$$
This is precisely the bound in
\eqref{eq:finite-main}, completing the proof.
\end{proof}

\section{Counting large central values}\label{sec:counting-large-values}

We apply Lemma~\ref{lem:short}, Proposition~\ref{prop:phase}, and Theorem~\ref{thm:finite} to estimate large
central values in the critical range. Throughout this section  $k\in\{3,19\}$. 
The integer $k$ is the energy parameter in Theorem~\ref{thm:finite}.
We take $k=3$ for the sixth moment, while
$k=19$ is used to match the critical-range
estimate with the high-value estimate of
Proposition~\ref{prop:tail} for the twelfth moment.

\begin{proposition}[Critical-range estimates]\label{prop:critical-range}
For every $\varepsilon>0$, the choice $k=3$ yields
\begin{equation}\label{eq:critical-claims}
R(q,q^\sigma)
\ll_{p,\varepsilon}
q^{(28-155\sigma)/18+\varepsilon},
\qquad
\frac{10}{83}\le\sigma\le\frac{13}{100},
\end{equation}
which will be used to improve the sixth moment bound.
Similarly, the choice $k=19$ gives
\begin{equation}
R(q,q^\sigma)
\ll_{p,\varepsilon}
q^{(1404-8315\sigma)/738+\varepsilon},
\qquad
\frac18\le\sigma\le\sigma_*,
\end{equation}
where
$\sigma_*=34794/261961$.
This estimate treats the critical range for the
twelfth moment.
\end{proposition}

\begin{proof}
Let $V=q^\sigma$ and $q_1$ be the largest $p$-power not
exceeding $V^2q^{-\varepsilon}$.
For each dyadic $F_0\ge1$, consider the
$H_{q_1}$-cosets containing between $F_0$
and $2F_0$ characters with
$|L(1/2,\chi)|\ge V$.
On each such coset,
$$
S_2(\chi;q_1)\ge F_0V^2,\qquad q_1\le V^2q^{-\varepsilon}.
$$
We apply Lemma~\ref{lem:short} with a sufficiently small
auxiliary parameter $0<\eta\ll\varepsilon$.
Since $q_1\le V^2q^{-\varepsilon}$, the
diagonal contribution satisfies
$$
q^\eta q_1\le V^2q^{\eta-\varepsilon}
=o(F_0V^2),
$$
uniformly for $F_0\ge1$, and can therefore
be absorbed. After absorbing the
$O(\log q)$ choices and normalizing the
coefficients, we obtain
\begin{equation}\label{eq:fixed-vector-large}
\left|\sum_{m\in I}a_mK_A^\pm(m;Q_h)\right|
\ge F_0Q_h^{1/2}q^{-\varepsilon},\qquad
|a_m|\le1,\quad |I|\ll q^\varepsilon N_h,
\end{equation}
where $D_h=q_1h$, $Q_h=q/D_h$, and $N_h=q/D_h^2$.
The coefficient system is independent of the center character;
its divisor-bound loss and the $O((\log q)^2)$ selections are absorbed in $q^\varepsilon$.

Write $h=q^u$ and $D_h=q^{2t}$. Then
\begin{equation}\label{eq:level-restriction}
t=\sigma+\frac{u}{2}-\varepsilon
+O(n^{-1}),\qquad
Q_h=q^{1-2t},\qquad N_h=q^{1-4t}.
\end{equation}
Since the left side of \eqref{eq:fixed-vector-large} is
$\ll q^\varepsilon N_h$, every nonempty layer satisfies
$N_h\gg Q_h^{1/2}q^{-\varepsilon}$, whence
$$
\sigma-\varepsilon\le t\le\frac16+ \varepsilon.
$$

We choose the modulus for Proposition~\ref{prop:phase}
and Theorem~\ref{thm:finite} by setting
\begin{equation}\label{eq:five-functions}
\begin{aligned}
\rho=\rho(t)
&=\min\left\{
\frac{4(k-1)+4(2k^2-6k+5)t}{(k-1)(2k+3)},
\frac43(1-4t)\right\},\\
H&=\max\{0,\rho-4t\}.
\end{aligned}
\end{equation}
Take
$$
R_h=p^{\lfloor(\rho-\varepsilon)n\rfloor},
\qquad q_{2,h}=D_hR_h.
$$
By Proposition~\ref{prop:phase},
$$
\eta=-\frac{D_h^2}{24B_0^2},\qquad
v_p(\eta)=2v_p(D_h)-\mathbf1_{p=3},
\qquad
\Upsilon_h=\frac{R_h}{(R_h,\eta)}
=q^{H+O(\varepsilon)}.
$$

For $t\in[10/83,1/6]$ if $k=3$, or
$t\in[1/8,1/6]$ if $k=19$, we have
\begin{equation}\label{eq:parameter-conditions}
\begin{gathered}
10t>1,\qquad
1-4t<\rho\le\frac43(1-4t),\\
2\rho<1+2t,\qquad
\rho<8t,\qquad
2t+\rho<1,\\
\rho+\min\{4t,\rho\}>2(1-4t).
\end{gathered}
\end{equation}
Indeed, checking the endpoints and the relevant
breakpoints $4/29$, $228/1597$, and $1/7$
shows that the least strict margin is $2/83$
for $k=3$ and $19/1476$ for $k=19$.
Since all strict inequalities in \eqref{eq:parameter-conditions} hold
with a uniform positive margin, they remain
valid in a sufficiently small neighborhood
of the stated intervals. As every nonempty
conductor-descent level satisfies
$\sigma-\varepsilon\le t\le 1/6+ \varepsilon$,
these inequalities continue to hold by
taking $\varepsilon>0$ sufficiently small.
The non-strict upper bound
$\rho\le \frac43(1-4t)$ follows directly
from \eqref{eq:five-functions}.

Since $|I|\ll q^\varepsilon N_h$, by slightly
enlarging the harmless $q^\varepsilon$ factor,
we may assume $|I|\le q^\varepsilon N_h$.
The strict margins in \eqref{eq:parameter-conditions} absorb this
adjustment and the prime-power rounding.
Thus Theorem~\ref{thm:finite} applies with length
parameter $q^\varepsilon N_h$, and
\begin{equation}\label{eq:actual-conditions}
\begin{gathered}
2\le q^\varepsilon N_h\le R_h
\le(q^\varepsilon N_h)^{4/3},\\
R_h(R_h,\eta)\ge q^{2\varepsilon}N_h^2,
\qquad D_hR_h<q,\\
5v_p(D_h)\ge n+\mathbf1_{p=5},\\
2v_p(R_h)\le n+v_p(D_h)-\mathbf1_{p=3},\\
v_p(R_h)\le4v_p(D_h)-2\mathbf1_{p=3}.
\end{gathered}
\end{equation}
Hence Proposition~\ref{prop:phase} and
Theorem~\ref{thm:finite} are applicable.

Fix an $H_{q_{2,h}}$-coset and a coefficient vector
as in \eqref{eq:fixed-vector-large}.
By Proposition~\ref{prop:phase}, absorbing the
$w$-independent phase and Gauss factor
into $a_m$ and replacing $\pm 2cw$ by $w$,
each selected short coset determines a
parameter $w$ modulo $R_h$ satisfying
$$
\left|\sum_{m\in I}a_m
e\left(\frac{w(r(m)+\eta r(m)^3)}{R_h}\right)
\right|\ge Z,
\qquad Z=F_0Q_h^{1/2}q^{-\varepsilon}.
$$
The map from short $H_{q_1}$-cosets to the
phase parameter $w\bmod R_h$ has fibers of
size $\phi(D_h)=h\phi(q_1)$. Hence, each fiber
consists of exactly $h$ short $H_{q_1}$-cosets. 
Therefore, if at most $|\mathcal W|$
parameters satisfy the above inequality, at most
$h|\mathcal W|$ short cosets contribute, and each containing fewer
than $2F_0$ large-value characters.
Since there are $O_p(q/q_{2,h})$ large cosets,
the total contribution from this dyadic class is
\begin{equation}\label{eq:fiber-count}
\ll_p \frac{q}{q_{2,h}}F_0h|\mathcal W|
=\frac{qF_0}{q_1R_h}\,|\mathcal W|,
\end{equation}
where $q_{2,h}=D_hR_h=q_1hR_h$.

By Theorem~\ref{thm:finite}, applied with length
$q^\varepsilon N_h$ and the above $Z$, we obtain
$$
\begin{aligned}
|\mathcal W|\ll q^\varepsilon\Bigg(
&\frac{N_h^2}{Z^2}
+\frac{N_h^{3-2/k}R_h^{(k+1)/(2k)}
\Upsilon_h^{1/(4k)}}{Z^{4-2/k}}\\
&+\frac{N_h^{k+1/2}R_h^{(k+1)/2}
\Upsilon_h^{1/4}}{Z^{2k+1}}
\Bigg).
\end{aligned}
$$
Since
$$
\begin{gathered}
N_h=q^{1-4t},\qquad Q_h=q^{1-2t},
\qquad q_1=q^{2t-u},\\
R_h=q^{\rho+O(\varepsilon)},\qquad
\Upsilon_h=q^{H+O(\varepsilon)},\qquad
Z=F_0Q_h^{1/2}q^{-\varepsilon},
\end{gathered}
$$
substituting into \eqref{eq:fiber-count} yields
\begin{equation}\label{eq:three-contributions}
\frac{qF_0}{q_1R_h}|\mathcal W|
\ll q^\varepsilon\left(
F_0^{-1}q^{A(t)+u}
+F_0^{-3+2/k}q^{B_k(t)+u}
+F_0^{-2k}q^{C_k(t)+u}
\right),
\end{equation}
where
\begin{equation}\label{eq:three-exponents}
\begin{aligned}
A(t)&=2-8t-\rho,\\
B_k(t)&=2-\frac1k-\left(10-\frac6k\right)t
-\frac{k-1}{2k}\rho+\frac H{4k},\\
C_k(t)&=1-(2k+3)t+\frac{k-1}{2}\rho
+\frac H4.
\end{aligned}
\end{equation}
Since $F_0\ge1$, all three powers of $F_0$ in
\eqref{eq:three-contributions} are at most one.
For the purpose of balancing, we regard
$B_k$ and $C_k$ in \eqref{eq:three-exponents} as functions
of the independent variables $t$ and
$\rho$, with $H=\max\{0,\rho-4t\}$.
When the first term in \eqref{eq:five-functions}
determines $\rho$,
we have $H=\rho-4t$ and
$$
\frac{\partial B_k}{\partial\rho}
=-\frac{2k-3}{4k}<0,\qquad
\frac{\partial C_k}{\partial\rho}
=\frac{2k-1}{4}>0.
$$
Thus the chosen $\rho$ balances $B_k$ and $C_k$,
giving
\begin{equation}\label{eq:balanced-exponent}
\begin{gathered}
B_k(t)=C_k(t)=
\frac{4k^2-2k-2-(24k^2-18k-7)t}
{(k-1)(2k+3)},\\
B_k(t)-A(t)=\frac{t}{k-1}>0.
\end{gathered}
\end{equation}
For $k=3,19$, after substituting
$\rho=\rho(t)$, differentiation on each
affine branch of \eqref{eq:five-functions} gives
\begin{equation}\label{eq:slope-bound}
\begin{aligned}
\max\{A'(t),B_k'(t),C_k'(t)\}
\le \max\bigg\{&
-\frac{24k^2-18k-7}{(k-1)(2k+3)},
-\frac83,\\
&-\frac{22}{3}+\frac{10}{3k},
-\frac{14k+1}{3}
\bigg\}<-2.
\end{aligned}
\end{equation}
Writing $h=q^u$, we have
$2t=\log_q q_1+u$, and hence
$$
\frac{d}{du}\bigl(G(t(u))+u\bigr)
=1+\frac12G'(t)<0,
\qquad G\in\{A,B_k,C_k\}.
$$
By continuity at the breakpoints, all three
exponents in \eqref{eq:three-contributions}
are maximized at $h=1$.
For $h=1$, our choice of $q_1$ gives
$t=\sigma-\varepsilon+O(n^{-1})$.
Moreover,
$$
\sigma\le\frac{13}{100}<\frac4{29}
\quad(k=3),\qquad
\sigma\le\sigma_*<\frac{228}{1597}
\quad(k=19).
$$
Thus $\rho$ lies on the balanced branch,
where $B_k=C_k>A$. Combining
\eqref{eq:three-contributions} and
\eqref{eq:balanced-exponent}, and summing over
the $O((\log q)^C)$ admissible choices, we obtain
\begin{equation}\label{eq:critical-intermediate}
R(q,q^\sigma)\ll_{p,k,\varepsilon}
q^{\frac{4k^2-2k-2-(24k^2-18k-7)\sigma}
{(k-1)(2k+3)}+\varepsilon}.
\end{equation}
Taking $k=3$ and $k=19$ yields, respectively,
the two bounds in Proposition \ref{prop:critical-range}.
\end{proof}

\section{Correlations and the high-value range}
\label{sec:high-values}
Let $p$ be a fixed odd prime and let $q=p^n$. Throughout this section,
$p\le q_1\mid q_2<q$ are powers of $p$ with $q_1\ge q^{1/5}$. For the
active powers $h=p^j$ arising from Lemma~3.1, we recall
$$
 D_h=q_1h,\qquad Q_h=\frac{q}{q_1h},\qquad
 N_h=\frac{q}{q_1^2h^2}.
$$
In particular, $h\mid Q_0$, where $Q_0=q/q_1$.

\subsection{An incomplete-correlation estimate}
\label{subsec:incomplete-correlation}

Let $A,A'$ be the Postnikov parameters of two primitive characters in
the same $H_{q_2}$-coset. Thus
$A\equiv A'\pmod{q/q_2}$, and in particular
$A\equiv A'\pmod p$. Define
\begin{equation}\label{eq:6.1}
\begin{gathered}
 Q_0=\frac{q}{q_1},\qquad N_0=\frac{q}{q_1^2},\qquad
 R_0=\frac{q_2}{q_1},\\
 P_0=\frac{Q_0}{(Q_0,A-A')},\qquad
 P_h=\frac{Q_h}{(Q_h,A-A')}
     =\max\left\{1,\frac{P_0}{h}\right\}.
\end{gathered}
\end{equation}
Indeed, if $h=p^j$ and $v=v_p(A-A')$, with $v_p(0)=+\infty$, then
$$
 P_0=p^{\max\{n-v_p(q_1)-v,0\}},\qquad
 P_h=p^{\max\{n-v_p(q_1)-j-v,0\}}.
$$
Let $\delta_q=(q/p,A-A')$. The effective modulus in
\cite[Proposition~2]{MW} is $Q_h/(Q_h,\delta_q)$; since
$Q_h\mid q/p$, this is precisely $P_h$. That proposition shows
that, for $P_h\ge p^2$, the same-sign product
$$
 G_{A,A'}^\pm(m;Q_h)
 :=K_A^\pm(m;Q_h)\overline{K_{A'}^\pm(m;Q_h)}
$$
is $P_h$-periodic. Moreover, its additive Fourier transform vanishes
unless the additive frequency is a $p$-adic unit. Taking the zero
frequency gives
$$
 \sum_{\substack{m\bmod P_h\\p\nmid m}}
 K_A^\pm(m;Q_h)\overline{K_{A'}^\pm(m;Q_h)}=0
 \qquad(P_h\ge p^2).
$$

\begin{lemma}\label{lem:corr}
Fix a $p$-adic exponent datum $\mathfrak d$ arising from the
$A$- and $B$-processes of \cite{Mil}, with exponent pair
$(\kappa,\lambda)$ satisfying
$0\le\kappa\le\frac12\le\lambda\le1$. There exists
$P_{\min}=P_{\min}(p,\mathfrak d)>p^2$ such that, if
$P_h>P_{\min}$ and $J\subset\mathbb Z$ is an interval containing
at most $M\ge1$ integers, then
\begin{equation}\label{eq:6.2}
 \left|\sum_{\substack{m\in J\\p\nmid m}}
 K_A^\pm(m;Q_h)\overline{K_{A'}^\pm(m;Q_h)}\right|
 \ll_{p,\mathfrak d,\varepsilon}
 q^\varepsilon P_h^\kappa M^{\lambda-\kappa}.
\end{equation}
For $P_h\le P_{\min}$, the trivial bound $O_p(M)$ holds.
\end{lemma}

\begin{proof}

Write $D=D_h$ and $Q=Q_h$. By the stationary-phase evaluation
in (2.6), both kernels are supported on the same admissible residue
classes, since $A\equiv A'\pmod p$. These classes are characterized by
$\left(\frac{mA^{-1}}p\right)=1.$
On all other unit classes both branches vanish. There are
$(p-1)/2=O_p(1)$ admissible classes, so we fix one of them and choose
compatible square-root branches. Write
$$
 A=cB^2,\qquad A'=c(B')^2,\qquad
 B\equiv B'\pmod p,
 \qquad c,B,B'\in\mathbb Z_p^\times.
$$
As $P_h>1$, we have $v_p(A-A')<v_p(Q)$. Write
$$
 g=p^{v_p(A-A')},\qquad Q=gP_h.
$$
Since $B+B'\in\mathbb Z_p^\times$, we have
$$
 A-A'=c(B-B')(B+B'),\qquad
 v_p(B-B')=v_p(g),
$$
so $(B-B')/g$ is a $p$-adic unit.

The stationary point $x_\pm=x_\pm(m)$ satisfies
$$
 Ax_\pm^2=m(1+Dx_\pm),\qquad
 \Phi_{A,D}^\pm(m)
 =\frac AD\log_p(1+Dx_\pm)+\frac m{x_\pm}.
$$
Differentiating the stationary value and using the critical-point
equation gives
\begin{equation}\label{eq:6.3}
 (\Phi_{A,D}^\pm)'(m)
 =\frac1{x_\pm(m)}
 =-\frac D2\pm\sqrt{\frac Am+\frac{D^2}{4}}.
\end{equation}
Let
$$
\mathscr P^\pm(m)
=
\Phi^\pm_{A,D}(m)-\Phi^\pm_{A',D}(m),
$$
where the same sign is used in both terms.
Expanding the square roots
on the fixed local branch yields
\begin{equation}\label{eq:6.4}
\begin{gathered}
 \frac{(\mathscr P^\pm)'(m)}g
 =\Omega_0\sqrt{\frac cm}
  \left(1+\sum_{j\ge1}D^{2j}b_jm^j\right),
 \qquad \Omega_0=\pm\frac{B-B'}g\in\mathbb Z_p^\times,\\
 b_j=\binom{1/2}{j}(4c)^{-j}
 \frac{B^{1-2j}-(B')^{1-2j}}{B-B'}.
\end{gathered}
\end{equation}
For every integer $a\ge1$, the finite geometric-sum identity gives
\begin{equation}\label{eq:6.5}
 \frac{B^{-a}-(B')^{-a}}{B-B'}
 =-\sum_{i=0}^{a-1}B^{-i-1}(B')^{-(a-i)}\in\mathbb Z_p.
\end{equation}
Since $p$ is odd, $\binom{1/2}{j}\in\mathbb Z_p$; hence
$b_j\in\mathbb Z_p$. The series in \eqref{eq:6.4} converges on the chosen
unit branch, because $p\mid D$.

Choose an integer $c_0\ge1$, to be fixed in terms of $p$ and
$\mathfrak d$, and write
$$
 m=m_0+p^{c_0}t,\qquad p\nmid m_0,
$$
on each class modulo $p^{c_0}$. Define
\begin{equation}\label{eq:6.6}
\begin{gathered}
 f_{m_0}(t)=\frac{\mathscr P^\pm(m_0+p^{c_0}t)-\mathscr P^\pm(m_0)}g,\\
 e\left(\frac{\mathscr P^\pm(m_0+p^{c_0}t)}Q\right)
 =e\left(\frac{\mathscr P^\pm(m_0)}Q\right)
  e\left(\frac{f_{m_0}(t)}{P_h}\right).
\end{gathered}
\end{equation}
The binomial expansions in \eqref{eq:6.4} give
$$
 \frac{(\mathscr P^\pm)'(m_0+p^{c_0}t)}g
 =\sum_{j\ge0}a_jt^j,
 \qquad v_p(a_j)\ge c_0j.
$$
Since $f_{m_0}(0)=0$,
$$
 f_{m_0}(t)
 =\sum_{j\ge0}\frac{p^{c_0}a_j}{j+1}t^{j+1}.
$$
Moreover,
$$
 v_p\left(\frac{p^{c_0}a_j}{j+1}\right)
 \ge c_0(j+1)-v_p(j+1)
 \ge 1,
$$
and the lower bound tends to $+\infty$ as $j\to\infty$.
Thus, $f_{m_0}\in\mathbb Z_p[[t]]$ and converges on $\mathbb Z_p$.
Let
$$
 G(z)=\sum_{j\ge1}D^{2j-2}b_jz^j,\qquad
 \omega=m_0^{-1},\qquad
 \omega'=\Omega_0\sqrt{c/m_0},\qquad
 u=2v_p(D).
$$
Here $c/m_0$ is a square on the chosen admissible class, so
$\omega'\in\mathbb Z_p^\times$. From \eqref{eq:6.4},
\begin{equation}\label{eq:6.7}
\begin{gathered}
 f_{m_0}'(t)
 =p^{c_0}\omega'(1+p^{c_0}\omega t)^{-1/2}
  +p^{c_0+u}g_1(t),\\
 g_1(t)=\omega'(1+p^{c_0}\omega t)^{-1/2}
  G(m_0+p^{c_0}t)
 =\sum_{j\ge0}g_{1,j}t^j,
 \qquad v_p(g_{1,j})\ge c_0j.
\end{gathered}
\end{equation}
Consequently, $g_1\in I_0[c_0](\mathbb Z_p)$ and
$
 f_{m_0}\in
 \mathcal F\left(c_0,\frac12,c_0,c_0,u,\omega,\omega'\right)
$
in the notation of \cite[Definition~1]{Mil}. In this notation,
$$
 y=\frac12,\quad w=\kappa_{\mathrm{fun}}
 =\lambda_{\mathrm{fun}}=c_0,\quad \gamma_0=0,
 \quad \iota(y)=\iota'(y)=0.
$$

We now check the hypotheses of \cite[Definition~2]{Mil}.
For odd $p$, $\iota'(2)=0$ and $\rho_p=1/(p-1)$; hence every
integer $c_0\ge1$ satisfies
$c_0\ge1+\iota'(2)$ and $c_0\in\rho_p\mathbb N$.
Choose $c_0$ sufficiently large that
$$
 c_0\ge\kappa_0(1/2,p),\qquad
 c_0\ge\lambda_0(1/2,p).
$$
The datum also gives
$n_0(1/2,p,c_0,c_0)>c_0$. Choose an integer
$C=C(p,\mathfrak d)$ sufficiently large that
$$
 C>2,\qquad
 C\ge c_0+n_0(1/2,p,c_0,c_0),\qquad
 \frac{2C}{5}\ge u_0(1/2,p,c_0,c_0),
$$
and let $P_{\min}=p^C$.
If $P_h=p^\nu>P_{\min}$, then
$\nu\ge c_0+n_0(1/2,p,c_0,c_0)\ge2c_0+1$. Moreover, since
$P_h\le Q_h\le q$ and $D=q_1h$,
$$
 u=2v_p(D)\ge\frac{2n}{5}>\frac{2C}{5}
 \ge u_0(1/2,p,c_0,c_0).
$$
Thus, all the hypotheses of Definition~2 are satisfied. Its permitted
interval length is
\begin{equation}\label{eq:6.8}
 B_{\mathrm{len}}\le p^{\nu-2c_0}.
\end{equation}
Applying \cite[Definition~2, (15)]{Mil}, we obtain
\begin{equation}\label{eq:6.9}
\begin{aligned}
 \left|\sum_{t\in J_t}
  e\left(\frac{f_{m_0}(t)}{P_h}\right)\right|
 &\ll_{\mathfrak d}
 p^{r_{\mathfrak d}}(p^{\nu-2c_0})^\kappa
 B_{\mathrm{len}}^{\lambda-\kappa}
 (\log p^{\nu-2c_0})^{\delta_{\mathfrak d}}\\
 &\ll_{p,\mathfrak d,\varepsilon}
 q^\varepsilon P_h^\kappa
 B_{\mathrm{len}}^{\lambda-\kappa}.
\end{aligned}
\end{equation}

The normalized stationary Gauss factors agree for the two
same-sign branches. Indeed, for
$\Theta_{A,D}(x,m)=A D^{-1}\log_p(1+Dx)+m/x$,
the critical-point equation gives
\begin{equation}\label{eq:6.10}
 \left.\frac12\frac{\partial^2\Theta_{A,D}}{\partial x^2}
 \right|_{x=x_*}
 =\frac{A(2+Dx_*)}{2x_*(1+Dx_*)^2}
 \equiv\frac{A}{x_*}\pmod p.
\end{equation}
Since $A\equiv A'\pmod p$ and the same-sign critical points agree
modulo $p$, the Gauss factors coincide and cancel in the conjugated
product. Thus \eqref{eq:6.9} applies without an additional amplitude.

If $M\le P_h=p^\nu$, each class modulo $p^{c_0}$ gives an interval of
$t$-length at most $1+M/p^{c_0}$. Since
$\nu\ge2c_0+1$, splitting these intervals into pieces satisfying
\eqref{eq:6.8} requires only $O_{p,\mathfrak d}(1)$ pieces. Applying \eqref{eq:6.9}
to each piece and using $\lambda\ge\kappa$ proves \eqref{eq:6.2}.
If $M>P_h$, remove the complete $P_h$-periods using the
complete-period identity above 
and apply the preceding estimate to
the remaining interval, whose length is less than $P_h$. Finally, for $P_h\le P_{\min}$, the normalized
stationary kernels satisfy $|K_A^\pm(m;Q_h)|\ll_p1$, so the
trivial bound is $O_p(M)$.
\end{proof}

\subsection{Aggregation of the short moments}
\label{subsec:aggregate}
For a character $\xi$ and an active level $h$, write
$$
 K_{\xi,h}^\pm(m):=K_{A_\xi}^\pm(m;Q_h),
$$
where $A_\xi$ is the Postnikov parameter of $\xi$.

\begin{proposition}\label{prop:6.2}
For every $\Psi\subseteq H_{q_2}$ and every primitive character
$\chi\pmod q$,
\begin{equation}\label{eq:6.11}
 \sum_{\psi\in\Psi}S_2(\chi\psi;q_1)
 \ll_{p,\mathfrak d,\varepsilon}q^\varepsilon
 \bigl(\bigl(q_1+q_1^{1/2}R_0^{\kappa/2}
 N_0^{(\lambda-\kappa)/2}\bigr)|\Psi|
 +q^{1/2}|\Psi|^{1/2}\bigr).
\end{equation}
\end{proposition}

\begin{proof}
By Lemma~\ref{lem:short}, summing over $\psi\in\Psi$
allows us to replace the maximum over
coefficient systems by their sum, at a cost
of $O(\log q)$ choices. We may therefore
fix one coefficient system, one stationary
sign, and one admissible class modulo $p$.
The coefficients are independent of the center character and satisfy
$$
 |a_h(m)|\ll_{p,\varepsilon}q^\varepsilon,\qquad
 \operatorname{supp}(a_h)\subseteq J_h,\qquad
 |J_h|\ll_{p,\varepsilon}q^\varepsilon N_h,
$$
where $J_h$ is an ordinary integer interval.
Throughout this subsection, sums over $m$
involving the stationary branches are
restricted to $p\nmid m$, and $a_h(m)$
is understood to vanish outside the
selected admissible residue class modulo $p$. 
Choose $|\gamma_{\psi,h}|\le1$ to linearize the corresponding
absolute values. Since $Q_h^{-1/2}=Q_0^{-1/2}h^{1/2}$, the
non-diagonal contribution is bounded by
$O(q_1Q_0^{-1/2}|T|)$, where
\begin{equation}\label{eq:6.12}
 T=\sum_hh^{1/2}\sum_m a_h(m)
    \sum_{\psi\in\Psi}\gamma_{\psi,h}
 K_{\chi\psi,h}^\pm(m),\qquad
 K_{\chi\psi,h}^\pm(m)=K_{A_{\chi\psi}}^\pm(m;Q_h).
\end{equation}
Cauchy's inequality gives
\begin{equation}\label{eq:6.13}
\begin{aligned}
 |T|^2\le{}&
 \left(\sum_hh^2\sum_m|a_h(m)|^2\right)
 \left(\sum_hh^{-1}\sum_m
 \left|\sum_{\psi\in\Psi}
 \gamma_{\psi,h}K_{\chi\psi,h}^\pm(m)\right|^2\right).
\end{aligned}
\end{equation}
Since $h^2N_h=N_0$ and there are $O(\log q)$ active levels,
$$
 \sum_hh^2\sum_m|a_h(m)|^2
 \ll_{p,\varepsilon}q^\varepsilon N_0.
$$
By nonnegativity, the second $m$-sum may first be extended to $J_h$.
Expanding the square only after this extension gives
$$
 \sum_{m\in\operatorname{supp}(a_h)}
 \left|\sum_{\psi\in\Psi}\gamma_{\psi,h}
 K_{\chi\psi,h}^\pm(m)\right|^2
 \quad\le\sum_{\psi,\psi'\in\Psi}
 \left|\sum_{\substack{m\in J_h\\p\nmid m}}
 K_{\chi\psi,h}^\pm(m)
 \overline{K_{\chi\psi',h}^\pm(m)}\right|.
$$
We omit levels with empty coefficient support.
For every remaining level, $1\le |J_h|
\ll q^\varepsilon N_h$, so Lemma~\ref{lem:corr}
applies with $M=|J_h|$.
Since $\lambda\ge\kappa$, the resulting
length losses are absorbed into $q^\varepsilon$.

For a fixed pair $(\psi,\psi')$, the quantities
$P_0$ and $P_h$ are understood to depend on
this pair. In the nondegenerate range
$P_h>P_{\min}$, Lemma~\ref{lem:corr} and $P_0\le R_0$
yield
\begin{equation}\label{eq:6.14}
\begin{aligned}
 \sum_{h:P_h>P_{\min}}h^{-1}P_h^\kappa N_h^{\lambda-\kappa}
 &\le R_0^\kappa N_0^{\lambda-\kappa}
 \sum_{h=p^j}h^{-(1+2\lambda-\kappa)}\\
 &\ll_p R_0^\kappa N_0^{\lambda-\kappa},\qquad
 P_0\le\frac{Q_0}{q/q_2}=R_0.
\end{aligned}
\end{equation}
For $P_h\le P_{\min}$, the trivial estimate and
$h\ge P_0/P_{\min}$ give
\begin{equation}\label{eq:6.15}
 \sum_{h:P_h\le P_{\min}}h^{-1}N_h
 \le N_0\sum_{\substack{h=p^j\\h\ge P_0/P_{\min}}}h^{-3}
 \ll_{p,\mathfrak d}N_0P_0^{-3}.
\end{equation}

For fixed $\psi$, the map $\psi'\mapsto A_{\chi\psi'}\bmod Q_0$
on $H_{q_2}$ has $R_0$ image points, each with $\phi(q_1)$
preimages. Moreover, for $P=p^j\le R_0$,
$$
 \begin{gathered}
 P_0(\psi,\psi')\le P
 \quad\Longleftrightarrow\quad
 A_{\chi\psi'}\equiv A_{\chi\psi}\pmod{Q_0/P},\\
 \#\{\psi'\in H_{q_2}:P_0(\psi,\psi')\le P\}
 \le\phi(q_1)P.
 \end{gathered}
$$
Consequently,
\begin{equation}\label{eq:6.16}
\begin{aligned}
 &\sum_{\psi,\psi'\in\Psi}
 \sum_{h:P_h\le P_{\min}}h^{-1}N_h\\
 &\quad\ll_{p,\mathfrak d,\varepsilon}q^\varepsilon
 N_0\sum_{\psi\in\Psi}
 \sum_{\substack{P=p^j\\P\le R_0}}
 P^{-3}\#\{\psi'\in\Psi:P_0(\psi,\psi')=P\}\\
 &\quad\ll_{p,\mathfrak d,\varepsilon}q^\varepsilon
 N_0\phi(q_1)|\Psi|\sum_{P=p^j}P^{-2}
 \ll_{p,\mathfrak d,\varepsilon}q^\varepsilon N_0q_1|\Psi|.
\end{aligned}
\end{equation}
Combining \eqref{eq:6.13}--\eqref{eq:6.16}, we have
\begin{equation}\label{eq:6.17}
 |T|^2\ll_{p,\mathfrak d,\varepsilon}q^\varepsilon N_0
 \left(R_0^\kappa N_0^{\lambda-\kappa}|\Psi|^2
       +N_0q_1|\Psi|\right).
\end{equation}
Finally, 
\begin{equation}\label{eq:6.18}
 q_1Q_0^{-1/2}N_0^{1/2}=q_1^{1/2},\qquad
 q_1Q_0^{-1/2}N_0q_1^{1/2}=q^{1/2}.
\end{equation}
Multiplying the square root of \eqref{eq:6.17} by $q_1Q_0^{-1/2}$ and
restoring the diagonal contribution $q_1|\Psi|$ proves \eqref{eq:6.11}.
\end{proof}

\subsection{A fixed exponent datum and the tail count}
\label{tail:count}

The choice $(\kappa,\lambda)=(\frac12,\frac12)$ in
\eqref{eq:6.11} recovers the coefficient
$q_1^{1/4}q_2^{1/4}$ of \cite[Proposition~3]{MW}.
We use the $p$-adic $A/B$-processes of
\cite[Theorems~4--5]{Mil}
\begin{equation}\label{eq:6.19}
 A(\kappa,\lambda)=
 \left(\frac{\kappa}{2(\kappa+1)},
       \frac{\kappa+\lambda+1}{2(\kappa+1)}\right),\qquad
 B(\kappa,\lambda)=
 \left(\lambda-\frac12,\kappa+\frac12\right).
\end{equation}
With composition from right to left,
\begin{equation}\label{eq:6.20}
 (\kappa_*,\lambda_*)
 =BA^7BABA^2BA^2B(0,1)
 =\left(\frac{2346}{4789},\frac{2400}{4789}\right),\qquad
 \lambda_*-\kappa_*=\frac{54}{4789}.
\end{equation}
The corresponding complete exponent datum is fixed henceforth.

\begin{proposition}\label{prop:tail}
We have
\begin{equation}\label{eq:6.22}
 R(q,V)\ll_{p,\varepsilon}
 q^{2-\delta_*+\varepsilon}V^{-12},\qquad
 q^{\sigma_*}\le V\le q^\beta,
\end{equation}
where
\begin{equation}\label{eq:6.21}
 \sigma_*=\frac{34794}{261961},\qquad
 \delta_*=\frac{51}{261961},\qquad
 \beta=\frac16-\frac1{1000}.
\end{equation}
\end{proposition}

\begin{proof}
Write $V=q^\sigma$, $\sigma_*\le\sigma\le\beta$.
For sufficiently small $\varepsilon>0$, let $q_1,q_2$ be the
largest powers of $p$ satisfying
\begin{equation}\label{eq:6.23}
 q_1\le V^2q^{-2\varepsilon},\qquad
 q_2\le q_1^3q^{\delta_*-\varepsilon}.
\end{equation}
The maximality of these powers gives
$$
 \begin{aligned}
 \log_q q_1&=2\sigma-2\varepsilon+O_p(1/\log q),\\
 \log_q q_2&=6\sigma+\delta_*-7\varepsilon+O_p(1/\log q),\\
 \log_q R_0&=4\sigma+\delta_*-5\varepsilon+O_p(1/\log q),\\
 \log_q N_0&=1-4\sigma+4\varepsilon+O_p(1/\log q).
 \end{aligned}
$$
Since $2\sigma_*>1/5$, $6\beta+\delta_*<1$, and
$\varepsilon<\delta_*$, we have $q_1\ge q^{1/5}$ and
$q_1\mid q_2<q$ for sufficiently large $q$.

By Proposition~\ref{prop:6.2}, it suffices to have
\begin{equation}\label{eq:6.24}
 R_0^{\kappa_*}N_0^{\lambda_*-\kappa_*}\le q_1.
\end{equation}
At the exponent level with $\varepsilon=0$, this is equivalent to
\begin{equation}\label{eq:6.25}
 \kappa_*\delta_*
 \le(2+4\lambda_*-8\kappa_*)\sigma
 -(\lambda_*-\kappa_*)
 \quad\Longleftrightarrow\quad
 \delta_*\le\frac{205\sigma-27}{1173},
 \qquad \frac{205\sigma_*-27}{1173}=\delta_*.
\end{equation}
The actual choices in \eqref{eq:6.23} give the strict margin
\begin{equation}\label{eq:6.26}
\begin{aligned}
 \log_q\left(
 \frac{R_0^{\kappa_*}N_0^{\lambda_*-\kappa_*}}{q_1}
 \right)
 &=-\frac{410}{4789}(\sigma-\sigma_*)
   -\frac{1936}{4789}\varepsilon
   +O_p(1/\log q)\\
 &<0.
\end{aligned}
\end{equation}
Indeed, $2+4\lambda_*-9\kappa_*=-1936/4789$.
Thus \eqref{eq:6.24} holds, and \eqref{eq:6.11} reduces to
\begin{equation}\label{eq:6.27}
 \sum_{\psi\in\Psi}S_2(\chi\psi;q_1)
 \ll_{p,\mathfrak d,\varepsilon}q^\varepsilon
 \left(q_1|\Psi|+q^{1/2}|\Psi|^{1/2}\right).
\end{equation}

Fix an $H_{q_2}$-coset $\chi H_{q_2}$ of
primitive characters, and let $R_2$ be
the number of characters in this coset
whose central values have absolute value
at least $V$.
Let $\Psi\subseteq H_{q_2}$ be the union
of those $H_{q_1}$-cosets whose translates
by $\chi$ meet this large-value set.
Then $\Psi H_{q_1}=\Psi$, and hence
\begin{equation}\label{eq:6.28}
\begin{gathered}
 |\Psi|\le\phi(q_1)R_2,\\
 \sum_{\psi\in\Psi}S_2(\chi\psi;q_1)
 =\phi(q_1)\sum_{\psi\in\Psi}|L(1/2,\chi\psi)|^2
 \ge\phi(q_1)R_2V^2.
\end{gathered}
\end{equation}
Using \eqref{eq:6.27}, the term proportional to
$q_1\phi(q_1)R_2$ is absorbed by
$q_1V^{-2}\le q^{-2\varepsilon}$. For $R_2>0$ this gives
$$
 R_2\ll_{p,\mathfrak d,\varepsilon}
 \frac{q^{1+\varepsilon}}{q_1V^4}.
$$
There are $O_p(q/q_2)$ primitive $H_{q_2}$-cosets, and
$$
 q_1q_2V^4
 \gg_p q_1^4q^{\delta_*-\varepsilon}V^4
 \gg_p V^{12}q^{\delta_*-9\varepsilon}.
$$
Consequently,
\begin{equation}\label{eq:6.29}
 R(q,V)
 \ll_{p,\mathfrak d,\varepsilon}
 \frac{q^{2+\varepsilon}}{q_1q_2V^4}
 \ll_{p,\varepsilon}
 q^{2-\delta_*+\varepsilon}V^{-12}.
\end{equation}
This proves \eqref{eq:6.22}.
\end{proof}

Finally, we record the upper cutoff for the central
values. By the pointwise estimate of
\cite{Mil}, for every
$\varepsilon>0$,
$$
|L(1/2,\chi)|
\ll_{p,\varepsilon}q^{27/164+\varepsilon}.
$$
Since $\beta>27/164$, choosing
$0<\varepsilon<\beta-27/164$ gives
$R(q,q^\beta)=0$ for sufficiently large $q$.

\section{Proof of the moment bounds}\label{sec:moments}

Let $\mathcal L_\chi=|L(1/2,\chi)|$, and recall that
$M_j(q)$ is defined in \eqref{intro:moments}. For $j>0$
and $1\le X<Y$, write
$$
[X,Y)=\coprod_{\substack{\ell\ge0\\2^\ell X<Y}}
[2^\ell X,\min\{2^{\ell+1}X,Y\}).
$$
On the $\ell$-th interval,
$\mathcal L_\chi^j\le2^j(2^\ell X)^j$. Hence
\begin{equation}\label{moments:dyadic-inequality}
\sum_{\substack{\chi\bmod q\\X\le\mathcal L_\chi<Y}}^{*}
\mathcal L_\chi^j
\le
2^j\sum_{\substack{\ell\ge0\\2^\ell X<Y}}
(2^\ell X)^jR(q,2^\ell X).
\end{equation}
By the pointwise estimate in \eqref{eq:inputs}, there are
$O(\log q)$ nonempty levels. 
If the last dyadic interval is truncated at an endpoint
of a counting range, only its left endpoint occurs in
$R(q,\cdot)$. 

\subsection{The sixth moment}\label{moments:sixth}

Let $V_0=q^{10/83}$ and $V_1=q^{13/100}$.
For $\mathcal L_\chi<V_0$, the fourth moment estimate
in \eqref{intro:classical} gives
\begin{equation}\label{eq:six-low}
\sum_{\substack{\chi\bmod q\\\mathcal L_\chi<V_0}}^{*}
\mathcal L_\chi^6
\le V_0^2M_4(q)
\ll_{p,\varepsilon}
 q^{1+20/83+\varepsilon}
=q^{103/83+\varepsilon}.
\end{equation}

For $V=q^\sigma$ with
$10/83\le\sigma\le13/100$, the first estimate of Proposition~\ref{prop:critical-range} gives
\begin{equation}\label{eq:six-middle}
V^6R(q,V)
\ll_{p,\varepsilon}
 q^{6\sigma+\frac{28-155\sigma}{18}+\varepsilon}=q^{\frac{28-47\sigma}{18}+\varepsilon}.
\end{equation}
Note that $\frac{28-47\sigma}{18}\le\frac{103}{83}$ when $10/83\le\sigma\le13/100$, with equality at $\sigma=10/83$.

For $\mathcal L_\chi\ge V_1$, we use only the previously
established twelfth moment estimate in \eqref{intro:classical} :
\begin{equation}\label{eq:six-high}
\sum_{\substack{\chi\bmod q\\\mathcal L_\chi\ge V_1}}^{*}
\mathcal L_\chi^6
\le V_1^{-6}M_{12}(q)
\ll_{p,\varepsilon}
 q^{2-78/100+\varepsilon}
=q^{61/50+\varepsilon}.
\end{equation}
Hence, after the dyadic summation, we obtain
$$
M_6(q)\ll_{p,\varepsilon}q^{103/83+\varepsilon}.
$$
In particular, this estimation uses only the
previously known twelfth moment bound in
\cite{MW}, not the improved bound proved in
the next subsection.

\subsection{The twelfth moment}\label{moments:twelfth}

We keep $V_0=q^{10/83}$ and the notation
$\sigma_*,\delta_*,\beta$ from \eqref{eq:6.21}.
For $\mathcal L_\chi<V_0$, the fourth moment estimate gives
\begin{equation}\label{eq:twelve-low}
\sum_{\substack{\chi\bmod q\\\mathcal L_\chi<V_0}}^{*}
\mathcal L_\chi^{12}
\le V_0^8M_4(q)
\ll_{p,\varepsilon}
 q^{1+80/83+\varepsilon}
=q^{163/83+\varepsilon}.
\end{equation}
For convenience, write
$E_3(\sigma)=(28-155\sigma)/18$ and
$E_{19}(\sigma)=(1404-8315\sigma)/738$,
as in Proposition~\ref{prop:critical-range}.
For $V=q^\sigma$, apply the two ranges of Proposition~\ref{prop:critical-range}:
\begin{equation}\label{eq:twelve-middle}
\begin{aligned}
12\sigma+E_3(\sigma)
&=\frac{28+61\sigma}{18}\le\frac{95}{48}=2-\frac{1}{48},
&&\frac{10}{83}\le\sigma\le\frac18,\\
12\sigma+E_{19}(\sigma)
&=\frac{1404+541\sigma}{738}
\le\frac{1404+541\sigma_*}{738}=2-\delta_* ,
&&\frac18\le\sigma\le\sigma_*.
\end{aligned}
\end{equation}
Both affine functions are increasing. Moreover,
$$
\begin{aligned}
\delta_*=\frac{51}{261961} <\frac{1}{48} <\frac{3}{83}.
\end{aligned}
$$
For $q^{\sigma_*}\le V\le q^\beta$,
Proposition~\ref{prop:tail} gives
$$
V^{12}R(q,V)
\ll_{p,\varepsilon}
q^\varepsilon V^{12}q^{2-\delta_*}V^{-12}
= q^{2-\delta_*+\varepsilon}.
$$
Since $\beta-27/164=127/123000>0$, the pointwise
estimate in \eqref{eq:inputs} implies
$$
\max_{\chi\bmod q}^{*}\mathcal L_\chi=o(q^\beta),
$$
and hence $R(q,q^\beta)=0$ for all sufficiently large $q$.

The transition identities are
\begin{equation}\label{eq:matching}
E_3\left(\frac{10}{83}\right)=1-4\times\frac{10}{83}=\frac{43}{83}, \qquad
\frac{1404+541\sigma_*}{738}
=2-\frac{205\sigma_*-27}{1173}=2-\delta_*.
\end{equation}
At $\sigma=1/8$ no matching identity is required:
the ranges $[10/83,13/100]$ and $[1/8,\sigma_*]$ overlap.
Thus \eqref{moments:dyadic-inequality},
\eqref{eq:twelve-low}, \eqref{eq:twelve-middle}, and
\eqref{eq:matching} cover every relevant dyadic level.
We conclude that
$$
M_{12}(q)
\ll_{p,\varepsilon}q^{2-\delta_*+\varepsilon}
=q^{2-51/261961+\varepsilon}.
$$
The finitely many excluded bounded depths are absorbed into the implied constants. This completes the proof of Theorem ~\ref{thm:main}.

\end{document}